\documentclass[12pt,twoside,leqno]{article}
\usepackage{amsmath}
\usepackage{amssymb}
\usepackage{amsxtra}
\usepackage{amscd}
\usepackage{amsthm}
\usepackage[mathscr]{eucal}

 \font\Russ=wncyr10   scaled\magstep 1

\def\Sha{\hbox{\Russ\char88}}

\theoremstyle{plain}
\newtheorem{thm}[subsection]{Theorem}
\newtheorem{prop}[subsection]{Proposition}
\newtheorem{cor}[subsection]{Corollary}
\newtheorem{lem}[subsection]{Lemma}
\newtheorem{conj}[subsection]{Conjecture}

\theoremstyle{definition}

\newtheorem{defn*}{Definition}

\newenvironment{pf}{\proof[\proofname]}{\endproof}

\begin{document}
\title{The $GL_2$ main conjecture for elliptic curves
without complex multiplication}
\author{J. Coates, T. Fukaya, K. Kato, R. Sujatha, O. Venjakob}
\maketitle

\newcommand\Cal{\mathcal}
\newcommand\define{\newcommand}

\define\gp{\mathrm{gp}}%
\define\fs{\mathrm{fs}}%
\define\an{\mathrm{an}}%
\define\mult{\mathrm{mult}}%
\define\Ker{\mathrm{Ker}\,}%
\define\Coker{\mathrm{Coker}\,}%
\define\Hom{\mathrm{Hom}\,}%
\define\Ext{\mathrm{Ext}\,}%
\define\rank{\mathrm{rank}\,}%
\define\gr{\mathrm{gr}}%
\define\Aut{\mathrm{Aut}\,}%
\define\cHom{\Cal Hom\,}%
\define\cExt{\Cal Ext\,}%
\define\cO{\Cal O}
\define\cS{\Cal S}
\define\cM{\Cal M}
\define\cG{\Cal G}
\define\cH{\Cal H}
\define\cE{\Cal E}
\define\cF{\Cal F}
\define\fm{\mathfrak m}

\newcommand{\N}{{\mathbb{N}}}
\newcommand{\Q}{{\mathbb{Q}}}
\newcommand{\Z}{{\mathbb{Z}}}
\newcommand{\R}{{\mathbb{R}}}
\newcommand{\C}{{\mathbb{C}}}
\newcommand{\F}{{\mathbb{F}}}
\newcommand{\bN}{{\mathbb{N}}}
\newcommand{\bQ}{{\mathbb{Q}}}
\newcommand{\bZ}{{\mathbb{Z}}}
\newcommand{\bR}{{\mathbb{R}}}
\newcommand{\bC}{{\mathbb{C}}}
\newcommand{\bF}{{\mathbb{F}}}
\newcommand{\bbQ}{{\bar \mathbb{Q}}}
\newcommand{\ol}[1]{\overline{#1}}
\newcommand{\too}{\longrightarrow}
\newcommand{\respect}{\rightsquigarrow}
\newcommand{\compatible}{\leftrightsquigarrow}
\newcommand{\upc}[1]{\overset {\lower 0.3ex \hbox{${\;}_{\circ}$}}{#1}}
\newcommand{\Gmlog}{\bG_{m, \log}}
\newcommand{\Gm}{\bG_m}
\newcommand{\eps}{\varepsilon}
\newcommand{\Spec}{\operatorname{Spec}}
\newcommand{\val}{{\mathrm{val}}}
\renewcommand{\b}{\langle\,,\,\rangle}
\newcommand{\n}{\operatorname{naive}}
\newcommand{\bs}{\operatorname{\backslash}}
\newcommand{\Gal}{\operatorname{\text{Gal}}}
\newcommand{\gal}{\operatorname{\text{Gal}({\bar \Q}/{\Bbb Q})}}
\newcommand{\galp}{\text{Gal}({\bar \Q}_p/{\Bbb Q}_p)}
\newcommand{\gall}{\text{Gal}({\bar \Q}_\ell/\Q_\ell)}
\newcommand{\GalQq}{\text{Gal}({\bar \Q}_q/{\Bbb Q}_q)}
\newcommand{\wep}{W({\bar \Q}_p/\Q_p)}
\newcommand{\wel}{W({\bar \Q}_\ell/\Q_\ell)}
\newcommand{\even}{\operatorname{even}}
\newcommand{\odd}{\operatorname{odd}}
\newcommand{\GL}{\operatorname{GL}}
\newcommand{\Lam}{{\Lambda}}
\newcommand{\La}{{\Lambda}}
\newcommand{\la}{{\lambda}}
\newcommand{\ga}{{\gamma}}
\newcommand{\om}{{\omega}}
\newcommand{\tLam}{{\tilde \Lambda}}
\newcommand{\tLa}{{\tilde \Lambda}}
\newcommand{\uL}{{{\hat L}^{\text{ur}}}}
\newcommand{\uQp}{{{\hat \Bbb Q}_p}^{\text{ur}}}
\newcommand{\sel}{\operatorname{\text{Sel}}}
\newcommand{\Om}{{\Omega}}
\newcommand{\OmG}{{\Omega(G)}}
\newcommand{\OmGp}{{\Omega'(G)}}
\newcommand{\OmGSp}{{\Omega(G)_{S'}}}
\newcommand{\OmH}{{\Omega(H)}}
\newcommand{\OmPi}{{\Omega(\Pi)}}
\newcommand{\OmGJ}{{\Omega(G/J)}}
\newcommand{\OmGJSJ}{{\Omega(G/J)_{S_J}}}
\newcommand{\LaG}{{\Lambda(G)}}
\newcommand{\LaOG}{{\Lambda_O(G)}}
\newcommand{\LaH}{{\Lambda(H)}}
\newcommand{\LaGH}{{\Lambda(G/H)}}
\newcommand{\LaHp}{{\Lambda(H')}}
\newcommand{\LaJ}{{\Lambda(J)}}
\newcommand{\LaPi}{{\Lambda(\Pi)}}
\newcommand{\LaGJ}{{\Lambda(G/J)}}
\newcommand{\LaGa}{{\Lambda(\Gamma)}}
\newcommand{\LaOGa}{{\Lambda_O(\Gamma)}}
\newcommand{\LaGS}{{\Lambda(G)_S}}
\newcommand{\LaGSa}{{\Lambda(G)_{S^*}}}
\newcommand{\LaU}{{\Lambda(U)}}
\newcommand{\vphi}{{\varphi}}
\newcommand{\vphiJ}{{\varphi_J}}
\newcommand{\psiJ}{{\psi_J}}
\newcommand{\RGa}{{R(\Gamma)}}
\newcommand{\ROGa}{{R_O(\Gamma)}}
\newcommand{\QGa}{{Q(\Gamma)}}
\newcommand{\QOGa}{{Q_O(\Gamma)}}
\newcommand{\MSaG}{{\frak {M}}_H(G)}
\newcommand{\Mf}{M_f}
\newcommand{\MO}{M_O}
\newcommand{\vare}{\varepsilon}
\newcommand{\mrho}{m_{\rho}}
\newcommand{\Delrho}{\Delta_{\rho}}
\newcommand{\Qcyc}{{{\mathbb{Q}}^{\rm {cyc}}}}
\newcommand{\Fin}{{F_{\infty}}}
\newcommand{\RSa}{R_{S^*}}
\newcommand{\cC}{{\cal {C}}}
\newcommand{\cD}{{\cal {D}}}
\newcommand{\fZ}{{\frak {Z}}}
\newcommand{\cA}{{\cal {A}}}
\newcommand{\cAG}{{\frak {A}(G)}}
\newcommand{\cR}{{\cal {R}}}
\newcommand{\cRG}{{\frak {R}(G)}}
\newcommand{\cLE}{{\cal {L}}_E}
\newcommand{\xiM}{{\xi_M}}
\newcommand{\Fcyc}{{F^{\rm {cyc}}}}
\newcommand{\Lcyc}{{L^{\rm {cyc}}}}
\newcommand{\XEFin}{{X(E/\Fin)}}
\newcommand{\YEFin}{{Y(E/\Fin)}}
\newcommand{\XEFcyc}{{X(E/\Fcyc)}}
\newcommand{\YEFcyc}{{Y(E/\Fcyc)}}
\newcommand{\XELcyc}{{X(E/\Lcyc)}}
\newcommand{\XEpFin}{{X(E'/\Fin)}}
\newcommand{\XEpFcyc}{{X(E'/\Fcyc)}}
\newcommand{\XEpLin}{{X(E'/\Lin)}}
\newcommand{\XEpLcyc}{{X(E'/\Lcyc)}}

\newcommand{\Ga}{\Gamma}
\newcommand{\Del}{\Delta}
\newcommand{\Sig}{\Sigma}
\newcommand{\sig}{\sigma}
\newcommand{\del}{\delta}
\newcommand{\bsig}{\bar{\sigma}}
\newcommand{\dla}{d_{\La}}
\newcommand{\can}{{\rm {can}}}
\newcommand{\dSla}{d_{S^{-1}\La}}
\newcommand{\tSla}{S^{-1}\La\otimes_{\La}}
\newcommand{\Sla}{S^{-1}\La}
\newcommand{\Isom}{{\rm {Isom}}}
\newcommand{\Mor}{{\rm {Mor}}}
\newcommand{\bg}{\bar{g}}
\newcommand{\id}{{\rm {id}}}
\newcommand{\cone}{{\rm {cone}}}
\newcommand{\al}{a}
\newcommand{\ChL}{{\cal{C}}(\La)}
\newcommand{\Image}{{\rm {Image}}}
\newcommand{\Dp}{{D_{\text{parf}}}}
\newcommand{\Zp}{{\Bbb Z}_p}
\newcommand{\Qp}{{\Bbb Q}_p}
\newcommand{\barQp}{{\bar{\Bbb Q}}_p}
\newcommand{\barQ}{{\bar{\Bbb Q}}}
\newcommand{\Zq}{{\Bbb Z}_q}
\newcommand{\Zl}{{\Bbb Z}_l}
\newcommand{\Ql}{{\Bbb Q}_l}
\newcommand{\Fp}{{\Bbb F}_p}
\newcommand{\Ak}{{\text {Ak}}}
\newcommand{\AkO}{{\rm {Ak}}_O}
\newcommand{\Akrho}{{\text {Ak}}_{\rho}}
\newcommand{\fp}{{\frak p}}
\newcommand{\ord}{{\text {ord}}}
\newcommand{\Phirho}{\Phi_{\rho}}
\newcommand{\tilPhirho}{\tilde{\Phi}_{\rho}}
\newcommand{\Phirhoetap}{\Phi'_{\rho \eta}}
\newcommand{\grho}{g_{\rho}}
\newcommand{\Sar}{\Phi'_{\rho}}
\newcommand{\tor}{{\rm {tor}}}
\newcommand{\Tor}{\operatorname{Tor}}
\newcommand{\twrho}{{\rm {tw}}_{\rho}}
\newcommand{\twhatrho}{{\rm {tw}}_{\hat{\rho}}}
\newcommand{\tweta}{{\rm {tw}}_{\eta}}
\newcommand{\twetan}{{\rm {tw}}_{\eta^n}}
\newcommand{\tw}{\operatorname{tw}}
\newcommand{\Ktor}{K^{\rm {tor}}}
\newcommand{\hatrho}{\hat{\rho}}
\newcommand{\Jac}{{\rm {Jac}}}
\newcommand{\GaF}{\Gamma_F}
\newcommand{\Epinf}{{E_{p^{\infty}}}}
\newcommand{\Frobq}{{\rm Frob}_q}
\newcommand{\Krho}{K_{\rho}}
\newcommand{\Vrho}{V_{\rho}}
\newcommand{\eprho}{e_p(\rho)}
\newcommand{\eqrho}{e_q(\rho)}
\newcommand{\Vrhola}{V_{\rho,\lambda}}
\newcommand{\Vl}{V_l}
\newcommand{\cV}{{\cal V}}

\section{Introduction}

The main conjectures of Iwasawa theory provide the only general method
known at present for studying the mysterious relationship between
purely arithmetic problems and the special values of complex
$L$-functions, typified by the conjecture of Birch and Swinnerton-Dyer and
its generalizations.  Our goal in the present paper is to develop
algebraic techniques which enable us to formulate a precise version of
such a main conjecture for motives over a large class of $p$-adic Lie
extensions of number fields. The methods which we develop in general were
inspired by the Heidelberg Habilitation thesis of one of us (Venjakob
\cite{Ve1}).

Let $G$ be a compact $p$-adic Lie group with no element of order $p$ and
write $\Lambda(G)$ for the Iwasawa algebra of $G$ (see \S 2). Let $M$ be a
finitely generated torsion $\Lambda(G)$-module and write $\chi(G,M)$ for
its $G$-Euler characteristic (see \S 3, (34)). How can we define a
characteristic element of $M$ and in particular relate
this to the Euler characteristic of $M$ and its twists? Let us
quickly recall how such characteristic elements are defined in classical
commutative Iwasawa theory when $G={\Bbb Z}_p^d$ for some integer $d \geq
1$. In this case, the structure theory for finitely generated torsion
$\Lambda(G)$-modules \cite{BCA}
shows that there exist a finite number of
non-zero elements $f_1,\cdots, f_r$ such that we have an exact sequence of 
$\Lambda(G)$-modules
$$
0 \to \oplus_{i=1}^r \Lambda(G)/\Lambda(G)f_i \to M \to D \to 0,
$$
where $D$ is a pseudo-null $\Lambda(G)$-module. We then define a
characteristic element of $M$  to be
$$
f_M:=\prod_{i=1}^r f_i(T)
$$ 
which is uniquely determined up to multiplication by a unit in
$\Lambda(G)$. The classical theory 
(see {\cite[V16]{30}}) shows that if
$H_0(G,M)$ is finite, then $\chi(G,M)$ and $\chi(G,D)$ are both finite,
and we have
$$
\chi(G,M)=\mid f_M(0)\mid^{-1}_p,~~~~\chi(G,D)=1
$$
where $f_M(0)$ denotes the image of $f_M$ under the augmentation map from
$\Lambda(G)$ to ${\Bbb Z}_p$. 

In the non-commutative case, assuming $G$ is $p$-valued, it is shown in
{\cite{31}} that there is an exact sequence 
$$
0 \to \oplus_{i=1}^r \Lambda(G)/L_i \to M/M_0 \to  D \to 0,
$$
where the $L_i$ are non-zero reflexive left ideals of $\Lambda(G)$,
$M_0$ is the maximal pseudo-null submodule of $M$, and $D$ is some
pseudo-null $\Lambda(G)$-module. The whole approach of the commutative
case to define characteristic elements now seems to break down
irretrievably. Firstly, it is no longer true (see {\cite{CSS}})  
that
$\chi(G,D)$ is finite implies that $\chi(G,D)=1$ for $D$ pseudo-null.
Secondly, it is also not true in general (see the appendix to 
\cite{32})
that a reflexive left ideal in $\Lambda(G)$ is always principal.  

The goal of this paper is to provide a way out of this dilemma via
localisation techniques for an important class of $G$. Namely, we assume
that $G$ has a closed normal subgroup $H$ such that $\Gamma=G/H \simeq
{\Bbb Z}_p$. For example, this is automatically true when $G$ is the
Galois group of a $p$-adic Lie extension of 
a number field $F$, 
which contains the cyclotomic ${\Bbb Z}_p$-extension of $F$.
We prove in \S 2 that the Iwasawa algebra $\Lambda(G)$ contains a
canonical Ore set $S^*$, enabling us to define the localised algebra
$\Lambda(G)_{S^*}$. Write ${\frak M}_H(G)$ for the category consisting of
all finitely generated $\Lambda(G)$-modules which are annihilated by $S^*$.
We prove that a finitely generated module $M$ belongs to ${\frak M}_H(G)$
if and only if $M/M(p)$ is finitely generated over $\Lambda(H)$ where
$M(p)$ denotes the $p$-primary submodule of $M$. We optimistically believe
that ${\frak M}_H(G)$ contains all the torsion $\Lambda(G)$-modules which
are of interest in arithmetic applications (for a precise statement see
Conjecture 5.1). We then exploit the well known localization sequence of
$K$-theory for the Ore set $S^*$ (see {\cite{Sw}})
to define a
characteristic element $\xi_M$ in $K_1(\Lambda(G)_{S^*})$ for any module 
$M$ in ${\frak M}_H(G)$. We also show that one can relate $\xi_M$ to the
Euler charcateristics of twists of $M$ by arbitrary continuous
representations of $G$, with values in $GL_n(O)$,  
where $O$ 
is the ring of integers of a finite extension of ${\Bbb Q}_p$. This uses
the Akashi series of $M$ which was introduced in {\cite{CSS}}.

The paper ends by formulating and briefly discussing the main conjecture
for an elliptic curve $E$ over ${\Bbb Q}$ over the field generated by the
coordinates of its $p$-power division points, where $p$ is a prime $\geq
5$ of good ordinary reduction for $E$. We also give some numerical
evidence in support of this main conjecture based on the remarkable
calculations of {\cite{27}}.
We have striven to keep the technical discussions to a minimum in the
present paper. A forthcoming paper (Fukaya and Kato, \cite{KF}) by two of us
will consider quite generally the Iwasawa theory of motives over arbitrary
$p$-adic Lie extensions of number fields and its connexion with the
Tamagawa number conjecture. In particular, that work applies to any
$p$-adic Lie extension with Galois group $G$, and does not need either of
the two basic hypotheses made in the present paper (namely that $G$ has no
element of order $p$, and that $G$ has a quotient isomorphic to $\Zp$).

\noindent{\it Acknowledgements:} This work was greatly assisted by the
generous hospitality provided by the Department of Mathematics of Kyoto
University for JC, the Harishchandra Research Institute, Allahabad for JC
and RS, and the Tata Institute of Fundamental Research, Mumbai for JC and
OV. 
Finally, we would like to thank Peter Schneider for pointing out to us the
alternative characterization of the Ore set $S$ given at the end of \S 2.\par

\section{The canonical Ore set.}

Let $G$ be a compact $p$-adic Lie group. We define
$$\La(G)=\lim_{\overset{\longleftarrow}{U}}\Zp[G/U], \quad
\Om(G)=\lim_{\overset{\longleftarrow}{U}}\Fp[G/U],   \eqno{(1)}\label{2.1}$$
where $U$ runs over all open normal subgroups of $G$.
By a module over these algebras, we shall always mean,
unless specified otherwise,
a left module.
It is well known that $\LaG$ and $\OmG$ are left and right Noetherian.
Moreover, if $\LaG$ has no zero divisors (e.g. if $G$ is
pro-$p$ and has no element of order $p$),
we shall always write $Q(G)$ for the skew field of fractions of $\LaG$.\par

We assume throughout this paper that $G$ has a closed normal subgroup $H$
such that
$$\Ga=G/H \overset{\simeq}{\to} \Zp.  \eqno{(2)} \label{2.2}$$
Inspired by the results in Venjakob \cite{Ve1}, our aim is to construct a
certain canonical
Ore set $S$ in $\LaG$.\par
\bigskip

{\bf {Definition.}}
Let $S$ be the set of all $f$ in $\LaG$ such that
$\LaG/\LaG f$ is a finitely generated $\LaH$-module.

\bigskip

Let $H'$ be any open subgroup of $H$. Then
$\LaH$ is a free left or right module of rank
$[H:H']$ over $\LaHp$.
Hence, if $M$ is a $\LaH$-module, $M$ will be finitely generated over $\LaH$ if
and only if it is finitely generated over $\LaHp$, and similarly if $M$ is
a right $\LaH$-module.
To exploit this observation, we will choose $H'=J$, where $J$ denotes any open
normal subgroup of $H$ which is pro-$p$.
For such a pro-$p$ $J$, $\LaJ$ is a local ring, whose maximal ideal is the
kernel of the
augmentation map from $\LaJ$ to $\Fp$, and
we can then use Nakayama's lemma for $\LaJ$-modules
(see the proof of Lemma \ref{l2.1}).
Let
$$\vphiJ:\LaG\to \LaGJ, \quad \psiJ:\LaG\to \OmGJ \eqno{(3)} \label{2.3}$$
be the natural surjections.
For $J$ as above, let us note that we can always find an open subgroup $J'$
of $J$
which is normal in $G$ because the set of $G$ conjugates of such a $J$ is
finite.

\begin{lem}
\label{l2.1}
Let $J$ be any pro-$p$ open subgroup of $H$, which is normal in $G$. Then
$\rm{(i)}$ $S$ is the set of all $f$ in $\LaG$ such that
$\LaGJ/\LaGJ \vphiJ(f)$ is a finitely generated $\Zp$-module,
$\rm{(ii)}$ $S$ is the set of all $f$ in $\LaG$ such that
$\OmGJ/\OmGJ\psi_J(f)$ is finite,
and $\rm{(iii)}$ $S$ is the set of all $f$ in $\La(G)$ such that right
multiplication by
$\psiJ(f)$ on $\OmGJ$ is injective.
\end{lem}

\begin{pf}
As remarked above, we can replace $H$ by $J$ in the definition of $S$. If
$f$ is any
element of $\LaG$, and if we put $M=\LaG/\LaG f$,
then we clearly have
$$(M)_J=\LaGJ/\LaGJ \vphiJ(f), \quad
M/\fm_J M=\OmGJ/\OmGJ\psiJ(f), \eqno{(4)} \label{2.4}$$
where $\fm_J$ denotes the maximal ideal of $\LaJ$.
Thus assertions (i) and (ii) are immediate from Nakayama's lemma.
To prove (iii), we use the important fact that we can always find a subgroup
$\Pi$ of $G/J$ satisfying
$$\Pi \overset{\simeq}{\to} \Zp, \quad
{\text {$\Pi$ is in the centre of $G/J$}}. \eqno{(5)} \label{2.5}$$
To establish this, let us write $\Ga'$ for some lifting of $\Ga$ to $G/J$.
Then we
can write $G/J$ as the semi-direct product of $H/J$ and $\Ga'$, where
$\Ga'$ acts on $H/J$
via conjugation. But, as $H/J$ is finite, an open subgroup of $\Ga'$ must
clearly act
trivially on $H/J$ by conjugation, and we take $\Pi$ to be this open subgroup.
    Now $\Om(\Pi)$ is isomorphic to the ring of formal power
series
$\Fp[[T]]$ in an indeterminate $T$ with coefficients in $\Fp$. Thus the
quotient field of
$\Om(\Pi)$, which we denote by
$R(\Pi)$, is a commutative field. Consider
$$V(G/J)=R(\Pi)\otimes_{\Om(\Pi)}\OmGJ=\OmGJ\otimes_{\Om(\Pi)}R(\Pi).
\eqno{(6)} \label{2.6}$$
It is a finite dimensional algebra over the commutative field $R(\Pi)$,
which lies in its
centre. If $f$ is any element of $\LaG$, let us write $\alpha_J(f)$ for the
element
of $V(G/J)$ defined by $\psiJ(f)$. By linear algebra, right multiplication by
$\alpha_J(f)$ is surjective if and only if it is injective.
Note also that, as $\OmPi$ is a discrete valuation ring with residue
field $\F_p$, a finitely generated
$\OmPi$-module is $\OmPi$-torsion if
and only if it is finite.
It follows easily that, for any $f$ in $\OmG$, $\OmGJ/\OmGJ \psiJ(f)$ is finite
if and only if right multiplication by $\psiJ(f)$ on $\OmGJ$ has finite kernel.
But $\OmGJ$ clearly has no finite $\OmGJ$-submodule, and so the equivalence
of (ii) and
(iii) is now clear.
\end{pf}

\begin{lem}
\label{l2.2}
Let $J$ be any pro-$p$ open subgroup of $H$, which is normal in $G$.
Then $\rm{(i)}$ $S$ is the set of all $f$ in $\La(G)$ such that
$\LaG/f\LaG$ is a finitely
generated
right $\LaH$-module; $\rm{(ii)}$ $S$ is the set of all $f$ in $\LaG$ such that
$\LaGJ/\vphiJ(f)\LaGJ$ is a finitely generated $\Zp$-module, $\rm{(iii)}$
$S$ is the set of
all
$f$ in $\LaG$ such that $\OmGJ/\psiJ(f)\OmGJ$ is finite, and $\rm{(iv)}$
$S$ is the set of
all
$f$ in
$\LaG$ such that the left multiplication by $\psiJ(f)$ on $\OmGJ$ is injective.
\end{lem}

\begin{pf} We first note that it suffices only to prove (iv). Indeed, once
(iv) is
established, we can simply reverse the arguments of the proof of Lemma
\ref{l2.1}, but
carrying them out for right modules, to deduce (i), (ii), and (iii). To
prove (iv), we
simply note the following. It was shown at the end of the proof of Lemma
\ref{l2.1} that $S$
consists of all
$f$ in $\LaG$ such that there exists $z$ in the algebra $V(G/J)$ with
$z.\alpha_J(f)=1$.
But, as multiplication on the left or on the right in $V(G/J)$ is
$R(\Pi)$-linear, and
$V(G/J)$ is finite dimensional over $R(\Pi)$, it follows from linear
algebra that
$z.\alpha_J(f)=1$ if and only if $\alpha_J(f).z=1$. This completes the
proof of lemma
\ref{l2.2}.
\end{pf}

\medskip

Let $M$ be a left or right $\LaG$-module. We say that $M$ is $S$-torsion
if, for each $x$ in
$M$, there exists $s$ in $S$ such that $s.x=0$ or $x.s=0$, according as the
action is on the
left or right.

\begin{prop} Let $M$ be a finitely generated left or right $\LaG$-module.
Then $M$ is
finitely generated over $\LaH$ if and only if $M$ is $S$-torsion.
\end{prop}

\begin{pf} We give the argument for left $\LaG$-modules. Suppose first that
$M$ is
$S$-torsion. Take a finite family $(x_i)_{1\leq i \leq r} $
of elements of $M$ which generates $M$ as a $\LaG$-module, and
choose $f_i$
in
$S$ such that
$f_i.x_i=0$ for each $i$.  Thus we get a
$\LaG$-surjection
$$\theta : \bigoplus_{i=1}^r \LaG/\LaG f_i \to M \eqno{(7)} \label{2.7}$$
by mapping $1$ in the $i$-th direct summand to $x_i$.
By the definition of $S$, each module occurring in the finite direct sum is
finitely
generated over $\LaH$, and hence $M$ is also. Conversely, assume that $M$
is finitely
generated over $\LaH$. We must show that every element $x$ of $M$ is
annihilated by an
element of $S$. As before, fix any pro-$p$ open subgroup $J$ of $H$, which
is normal in $G$.
Again, we choose a subgroup $\Pi$ of $G/J$ satisfying (5),
and we identify
$\La(\Pi)$ with $\Zp[[T]]$ by fixing a topological generator of $\Pi$. Now
$M$ is finitely
generated over $\LaJ$.
Let $\tau$ denote any lifting of $T$ to $\LaG$.
For each integer $n\geq 1$, define the $\LaJ$-submodule
$$U_n = \LaJ x + \LaJ \tau x + \cdots + \LaJ \tau^nx. \eqno{(8)} \label{2.8}$$
Since $\LaJ$ is Noetherian, it follows that there must exist an integer
$n\geq 1$, and $a_0,
\dots, a_{n-1}$ in $\La(J)$ such that
$$\tau^nx= (a_0+a_1\tau+\cdots+a_{n-1}\tau^{n-1})x.$$
Hence, if we define
$$s_n=\tau^n-a_{n-1}\tau^{n-1}-\cdots - a_0,$$
we have $s_n.x=0$. But clearly $\psi_J(s_n)$ is a non-zero element of
$\OmPi$, and hence,
by (iii) of lemma 2.1, $s_n$ belongs to $S$. This completes the proof of
Proposition 2.3.
\end{pf}

\begin{thm}
The set $S$ is multiplicatively closed, and is a left and right Ore set in
$\LaG$.
The elements of $S$ are non-zero divisors in $\LaG$.
\end{thm}

\begin{pf}
Take $s_1$ and $s_2$ in $S$. Then we have the exact sequence of $\LaG$-modules
$$0 \to W \to \LaG/\LaG s_1s_2 \to \LaG/\LaG s_2 \to 0, \eqno{(9)}
\label{2.9}$$
where $W=\LaG s_2/\LaG s_1s_2$. As $W$ is a homomorphic image of $\LaG/
\LaG s_1$,
it follows that the $\LaH$-module in the middle of (9) is finitely generated
over $\LaH$, and thus $s_1s_2$ belongs to $S$. Now take $f$ to be any
element of $S$,
and $x$ to be any element of $\LaG$.
By Proposition 2.3, the left module $\LaG/\LaG f$
and the right module $\LaG/f\LaG$ are both $S$-torsion.
Hence there exist $s$ and $s'$ in $S$ such that $sx \in \LaG f$ and
$xs' \in f \LaG$. This proves that $S$ is a left and right Ore set in $\LaG$.

To prove the final assertion of the theorem, take $f$ to be any element of $S$.
Now
$$\LaG = \lim_{\overset{\longleftarrow}{J}} \LaGJ,$$
where the projective limit is taken over all pro-$p$ open subgroups $J$ of $H$,
which are normal in $G$. Thus, fixing any such $J$,
it suffices to show that $\vphiJ(f)$ is not a zero divisor in $\LaGJ$.
As before, we choose a subgroup $\Pi$ of $G/J$ satisfying (5),
and let $Q(\Pi)$ denote the quotient field of $\La(\Pi)$.
Consider the algebra
$$W(G/J)=Q(\Pi)\otimes_{\La(\Pi)}\LaGJ=\LaGJ\otimes_{\La(\Pi)}Q(\Pi),
\eqno{(10)} \label{2.10}$$
which is finite dimensional over the commutative field $Q(\Pi)$, which lies
in its centre.
By virtue of (i) of Lemma \ref{l2.1} and (ii) of Lemma \ref{l2.2}, we see
that the image of
$\vphiJ(f)$ in $W(G/J)$ has a right and left inverse.
Hence $\vphiJ(f)$ cannot be a divisor of zero in $\LaGJ$, and
the proof of Theorem 2.4 is complete.
\end{pf}

We are grateful to P. Schneider (Lemma 2.5 and Proposition 2.6 below are
due to him) for kindly pointing out to us the following alternative
characterisation of the set $S$.\par

\begin{lem}
\label{l2.5}
For any two pro-p open subgroups $J \subseteq J'$ of $H$ which
are normal in $G$, the kernel of the natural map $\Omega(G/J) \to
\Omega(G/J')$ is a nilpotent ideal.
\end{lem}

\begin{pf}
We recall that, quite generally, the prime radical
${\cal N}(A)$ in
a Noetherian ring $A$ is nilpotent and contains any other nilpotent ideal
\cite[0.2.6, 2.3.7]{MR}. The assertion therefore is the claim that
$$
\Ker(\Omega(G/J) \to \Omega(G/J') \subseteq {\cal N}(\Omega(G/J)).
$$
The subgroup $\Delta:=J'/J$ is a finite normal $p$-subgroup of $G/J$. The
elements $\delta-1$ for $\delta \in \Delta$ generate the above left hand
side as a left ideal. Let $\Pi \subseteq G/J$ be a central subgroup as in
the proof of Lemma 2.1; it maps isomorphically onto a central subgroup in
$G/J'$. Hence we may view $\Omega(G/J) \to \Omega(G/J')$ as a map of
$\Omega(\Pi)$-algebras. Since
$$
{\cal N}(\Omega(G/J))=\Omega(G/J) \cap {\cal N}(V(G/J))
$$
it suffices to prove that
$$
\Ker(V(G/J) \to V(G/J')) \subseteq {\cal N}(V(G/J)).
$$
Again, the left hand side is generated, as a left ideal, by the elements 
$\delta -1$ for $\delta \in \Delta$. The right hand side is the
intersection of the  annihilator (prime) ideals of all the simple modules
over the finite dimensional $R(\Pi)$-algebra $V(G/J)$. Hence we have to 
show that $\Delta$ acts trivially on any simple $V(G/J)$-module $E$. Since
$\Delta$ is a $p$-group and $E$ is a vector space over a field of
characteristic $p$, the fixed vectors $E^{\Delta}$ under $\Delta$ in $E$
certainly are non-zero. Since $\Delta$ is normal in $G/J$, these fixed
vectors $E^{\Delta}$ form a $V(G/J)$-submodule of $E$. Hence we must have
$E^{\Delta}=E$.
\end{pf}

The lemma implies that the two sided ideal 
$$
{\cal N}:={\rm preimage}~{\rm of}~{\cal N}(\Omega(G/J))~{\rm in}~
\Lambda(G)
$$
is independent of the choice of a pro-$p$ open subgroup $J$ in $H$ which
is normal in $G$.

\begin{prop} 
\label{p2.6}
The set $S$ is equal to the set of all elements in
$\Lambda(G)$ which are regular modulo ${\cal N}$.
\end{prop}

\begin{pf}
Choose a pro-$p$ open subgroup $J \subseteq H$ which is
normal in $G$. Since $\Lambda(G)/{\cal N}
=\Omega(G/J)/{\cal N}(\Omega(G/J))$ 
an element
$f \in \Lambda(G)$ is regular modulo ${\cal N}$ if and only if $\psi_J(f)$
is regular modulo ${\cal N}(\Omega(G/J))$,
i.e.
becomes a unit in the finite dimensional algebra $V(G/J)/{\cal N}(V(G/J)).$ 
But units modulo nilpotent ideals are units. Hence $f$ is regular modulo
${\cal N}$
if and only if $\Psi_J(f)$ is a non-zero divisor in $\Omega(G/J)$. 
By Lemma 2.1(iii) and 2.2(iii), the latter is equivalent to $f$ belonging
to $S$.
\end{pf}

\section{Akashi series and Euler characteristics.}

Most of the results in this section are already established in \cite{Ve1},
but we reprove them here in a slightly different fashion.
As always, $G$ will denote a compact $p$-adic Lie group with a closed normal
subgroup $H$
such that $G/H=\Ga$ is isomorphic to $\Zp$.
We fix from now on a topological generator of $\Ga$, and identify $\LaGa$
with the formal
power series ring $\Zp[[T]]$ by mapping
this topological generator to $1+T$. We write
$\QGa$ for the fraction field of $\LaGa$.
Similarly, if $O$ denotes the ring of integers of some finite extension of
$\Qp$, we write
$\LaOGa$ (which we identify with $O[[T]]$)
for the $O$-Iwasawa
algebra of $\Ga$, and $\QOGa$ for the quotient field of $\LaOGa$.
Let $S$ be the Ore set in $\LaG$ which is defined in \S 2. Since $p \notin S$,
we shall also need to consider

\bigskip

{\bf {Definition.}}
$$S^*= \bigcup_{n \geq 0}p^nS.$$
\bigskip

\noindent
As $p$ lies in the centre of $\LaG$, $S^*$ is again a multiplicatively
closed left and right Ore set in $\LaG$,
all of whose elements are non-zero divisors.
We write $\LaGS$, $\LaGSa$ for the localizations of $\LaG$ at $S$ and
$S^*$, so that

$$\LaGSa=\LaGS[\frac{1}{p}]\eqno{(11)}$$
If $M$ is a $\LaG$-module, we write $M(p)$
for the submodule of $M$ consisting of all elements of finite order.
It is clear from Proposition 2.3 that $M$ will be $S^*$-torsion if and only
if $M/M(p)$
is finitely generated over $\LaH$.
We write $\MSaG$ for the category of all finitely generated $\LaG$-modules,
which are $S^*$-torsion. Note that in the special case in which 
$H=1$ and $G=\Ga$, 
$\MSaG$ is the category of all finitely generated torsion
$\LaG$-modules.
Both for motivation,
and because we shall need it later in this section, we prove the following
lemma
(see \cite{CSS}).

\begin{lem}
\label{l3.1}
For each $M$ in $\MSaG$, the homology groups $H_i(H,M)$ $(i \geq 0)$
are all finitely generated torsion $\LaGa$-modules. If $G$ has no element
of order $p$,
$H_i(H,M)=0$ for $i\geq d$, where $d$ is the dimension of $G$ as a $p$-adic Lie
group.
\end{lem}

\begin{pf}
We first observe that the $H_i(H,M)$ are
finitely generated $\LaGa$-modules, because
$$H_i(H,M)=\Tor_i^{\LaG}(\LaGH,M) \qquad (i \geq 0),$$
and the modules on the right are finitely
generated over $\LaGa$ since $M$ is finitely generated over $\LaG$.
Put $\Mf=M/M(p)$.
Now, for all $i \geq 0$, $H_i(H, M(p))$ is killed by $p^t$,
where $t$ is any integer $\geq 0$ such that $p^t. M(p)=0$.
Also, $H_i(H, \Mf)$ is a finitely generated $\Zp$-module since
$\Mf$ is finitely generated over $\LaH$.
It now follows from the long exact sequence of $H$-homology that
$H_i(H,M)$ is a torsion $\LaGa$-module for all $i \geq 0$.
The final assertion of the lemma is true because $H$ has
$p$-cohomological dimension $d-1$ when
$H$ has no elements of order $p$ \cite{Se2}.
\end{pf}

As above, let $O$ denote the ring of integers of some finite extension $L$
of $\Qp$,
and let us assume that we are given a continuous homomorphism
$$\rho:G \to GL_n(O), \eqno{(12)}$$
where $n$ is some integer $\geq 1$.
If $M$ is a finitely generated $\LaG$-module, put
$\MO=M\otimes_{\Zp}O$, and define
$$\twrho(M)=\MO\otimes_OO^n.  \eqno{(13)}$$
We endow $\twrho(M)$ with the diagonal action of $G$
i.e. if $\sig$ is in $G$, $\sig(m \otimes z)=(\sig m)\otimes (\sig z)$,
where it is understood that $G$ acts on $O^n$ on the left via
the homomorphism $\rho$.
By compactness, this left action of $G$ extends to an action of
the whole Iwasawa algebra $\LaG$.

\begin{lem}
\label{l3.2}
If $M \in \MSaG$, then, for all continuous representations
$\rho$ of the form $\rm{(12)}$, $\twrho(M) \in \MSaG$.
\end{lem}

\begin{pf}
Assume that $M \in \MSaG$.
Since $O$ is a free module of finite rank over $\Zp$,
it follows easily that $\MO \in \MSaG$.
Since $\twrho(\MO(p))$ is killed by the same power of $p$ which kills
$\MO(p)$, it suffices to prove that
$\twrho(R)$ is finitely generated over $\LaH$,
where $R=\MO/\MO(p)$. Now we have a surjection
$\LaH^m\to R$, which plainly induces a surjection of
$\LaH$-modules $\tw_{\rho_H}(\LaH^m)\to \twrho(R)$,
where $\rho_H$ denotes the restriction of $\rho$ to the subgroup $H$
of $G$. But it is well known \cite{Ve1} that
$\tw_{\rho_H}(\LaH^m)$ is again a free $\LaH$-module of finite rank.
This completes the proof of Lemma 3.2.
\end{pf}

If $R$ is a ring with unit element, we write $R^{\times}$
for the group of units of $R$, and $M_n(R)$ for
the ring of $n \times n$-matrices
with entries in $R$. By continuity,
the group homomorphism $\rho$ induces a ring homomorphism
$$\rho:\LaG\to M_n(O). \eqno{(14)}$$
A second ring homomorphism
$$\Phirho:\LaG \to M_n(\LaOGa), \eqno{(15)}$$
which we now define, will play an important role in our work.
If $\sig$ is in $G$, we write $\bsig$ for its image in
$\Ga=G/H$. We define a continuous group homomorphism
$$G \to (M_n(O) \otimes_{\Zp}\LaGa)^{\times} \eqno{(16)}$$
by mapping $\sig$ to $\rho(\sig) \otimes \bsig$.
Noting that
$$M_n(O)\otimes_{\Zp}\LaGa=M_n(\LaOGa)$$
because $O$ is a free $\Zp$-module of finite
rank, we obtain (15) by extending (16) to the whole of $\LaG$.
Note that $\Phirho(p)=pI_n$, where $I_n$ is the unit matrix.
Also it is easily seen that (14) is the composition of (15) with
the map from $M_n(\LaOGa)$ to $M_n(O)$ induced by
the augmentation map from $\LaOGa$ to $O$.

\begin{lem}
\label{l3.3}
The map $\rm{(15)}$ extends to a ring homomorphism, which we also
denote by $\Phirho$,
$$\Phirho:\LaGSa\to M_n(\QOGa). \eqno{(17)}$$
\end{lem}

\begin{pf}
We must show that, for all $s$ in $S^*$,
$\Phirho(s)$ is invertible in $M_n(\QOGa)$, or equivalently
that $\Phirho(s)$ has non-zero determinant.
Since this is clearly true for $s=p$, we can assume that $s$
is in $S$. Let $k$ denote the residue field of $O$, and let
$$\tilPhirho:\LaG\to M_n(k[[T]])$$
denote the composition of
$\Phirho$ with the canonical map from
$M_n(\LaOGa)$ to $M_n(k[[T]])$.
It suffices to show that
$\tilPhirho(s)$ has non-zero determinant. Note that
$\tilPhirho$ is induced by the map
$$G \to (M_n(k)\otimes_kk[[T]])^{\times} \eqno{(18)}$$
given by $\sig \mapsto \tilde{\rho}(\sig) \otimes \bsig$, where
$\tilde{\rho}(\sig)$ denotes the image of $\rho(\sig)$
in $M_n(k)$.
We now exploit the fact that $GL_n(k)$ is a finite group.
Thus we can find a pro-$p$ open subgroup $J$ of $H$,
which is normal in $G$, and which is contained in
$\Ker(\tilde{\rho})$.
Clearly $\tilPhirho$ can be factored through a map
$$\del:\Om(G/J) \to M_n(k[[T]]). \eqno{(19)}$$
Let $\psiJ$ be given by (3). Then we must show that
$\del(\psiJ(s))$ has non-zero determinant.
As in \S 2, we can find a central subgroup
$\Pi$ of $G/J$ with $\Pi \overset{\simeq}{\to} \Zp$.
Replacing $\Pi$ by an open subgroup if necessary,
we may also assume that $\Pi$ lies in the kernel of the homomorphism
from $G/J$ to $GL_n(k)$ induced by $\tilde{\rho}$.
Hence it is clear from (18) that
$\del(\OmPi)$ must be contained in $k[[T]].I_n$.
Let $R(\Pi)$ (resp. $Q(k[[T]])$) denote the quotient field of $\OmPi$
(resp. $k[[T]]$).
Since the natural map from $\Pi$ to $\Ga$ is injective,
$Q(k[[T]])$ can be viewed as a finite extension of $R(\Pi)$.
Hence (19) induces a ring homomorphism
$$\alpha:R(\Pi)\otimes_{\OmPi}\OmGJ \to M_n(Q(k[[T]])). \eqno{(20)}$$
But, as proven in \S 2, $\psiJ(s)$ is a unit in
the algebra on the left in (20).
Hence $\alpha(\psiJ(s))$ is invertible in
$M_n(Q(k[[T]])$, and so has non-zero determinant.
This completes the proof
of Lemma 3.3.
\end{pf}

If $R$ is a ring, we write $K_mR$ for the $m$-th
$K$-group of $R$
(we shall only need the cases $m=0,1$).
Clearly (14) induces a homomorphism
$$K_1(\LaG)\to K_1(M_n(O))=O^{\times}. \eqno{(21)}$$
This can be extended to a map from $K_1(\LaGSa)$ to
$L \cup \{\infty\}$, where $L$ is the fraction field of $O$,
in the following manner.
Firstly, (17) induces a homomorphism
$$\Sar:K_1(\LaGSa) \to K_1(M_n(\QOGa)))=\QOGa^{\times}.\eqno{(22)}$$
Let $\vphi:\LaOGa \to O$ be the augmentation map,
and write $\fp=\Ker(\vphi)$.
Of course, writing $\LaOGa_{\fp} \subset \QOGa$ for the
localization of $\LaOGa$ at $\fp$, it is clear that
$\vphi$ extends naturally to a homomorphism
$$\vphi:\LaOGa_{\fp}\to L.$$
Let $\xi$ be any element of $K_1(\LaGSa)$.
If $\Sar(\xi)$ belongs to $\LaOGa_{\fp}$, we define
$\xi(\rho)$ to be $\vphi(\Sar(\xi))$.
However, if $\Sar(\xi)$ does not belong to $\LaOGa_{\fp}$,
we define $\xi(\rho)=\infty$.
This gives us the desired extension of (21). \par

We now use a well known localization theorem in $K$-theory
to define the notion of a characteristic element for any
module $M$ in the category $\MSaG$.
To ensure that we can work with modules rather than complexes,
we assume for the rest of this section the following\par

\bigskip

{\bf {Hypothesis on $G$.}}
\qquad $G$ has no element of order $p$.

\bigskip

It is well known \cite{Bru} that
this implies that $\LaG$ has finite
global dimension equal to $d+1$, where $d$ is the dimension of $G$ as
a $p$-adic Lie group.
Since $\MSaG$ is the category of all finitely generated
$\LaG$-modules which are $S^*$-torsion, there is a connecting
homomorphism (see \cite{Sw})
$$\partial_G:K_1(\LaGSa)\to K_0(\MSaG), \eqno{(23)}$$
where $K_0(\MSaG)$ denotes the Grothendieck group of the category
$\MSaG$, such that we have, as part of a larger exact sequence
of localization, the exact sequence
$$\cdots \to K_1(\LaG)\to K_1(\LaGSa)\overset{\partial_G}{\to}
K_0(\MSaG)
\to K_0(\LaG)\to K_0(\LaGSa)\to 0.
\eqno{(24)}$$

\begin{prop}
\label{p3.4}
Assume that $G$ has no element of order $p$.
Then $\partial_{G}$ is surjective.
\end{prop}

Before giving the proof of Proposition 3.4,
we need the following preliminary lemma.\par

\begin{lem}
\label{l3.5}
Let $P$ be any pro-$p$ open normal subgroup of $G$.
Then the canonical map
$$K_0(\LaG) \to K_0(\Zp[G/P])$$
is injective.
\end{lem}

\begin{pf}
Let $\Del=G/P$, and let $I$ denote the kernel of
the natural map from $\LaG$ to $\Zp[\Del]$.
Let $M$ be any finitely generated projective
$\LaG$-module such that the class of $M/IM$
in $K_0(\Zp[\Del])$ is zero, i.e. we have a
$\Zp[\Del]$-isomorphism
$$\alpha:M/IM \oplus \Zp[\Del]^r \overset{\simeq}{\to}
\Zp[\Del]^r, $$
for some integer $r \geq 0$. Since $M$ is a projective
$\LaG$-module, we can find a $\LaG$-homomorphism
$$\beta:M\oplus \LaG^r \to \LaG^r \eqno{(25)}$$
which lifts $\alpha$. In particular, it is then
clear that $(\Coker(\beta))_P=0$, whence, as $P$
is pro-$p$, it follows from Nakayama's lemma that
$\Coker(\beta)=0$. Taking $P$-homology of the
exact sequence
$$0 \to \Ker(\beta) \to M\oplus \LaG^r
\to \LaG^r \to 0,$$
and noting that
$H_1(P, \LaG^r)=0$ because
$\LaG$ is a free $\La(P)$-module of finite rank,
we conclude that $(\Ker(\beta))_P=0$, and so also
$\Ker(\beta)=0$. Thus $\beta$ is an isomorphism, and the class
of $M$ in $K_0(\LaG)$ is zero. This completes the proof of
Lemma 3.5.
\end{pf}

We now prove Proposition 3.4.
For the proof, we fix a pro-$p$ open normal subgroup $P$ of $G$,
and put $\Del=G/P$. We write $\cV=\cV(\Del)$ for the
set of irreducible representations of the finite group $\Del$
over $\barQp$,
and we take $L$ to be some fixed finite extension of $\Qp$
such that all representations in $\cV$ can be realized over $L$.
Thus we have an isomorphism of rings
$$\gamma:L[\Del]\overset{\simeq}{\to}\prod_{\rho \in \cV}M_{n_{\rho}}(L),
\eqno{(26)}$$
where $n_{\rho}$ denotes the dimension of $\rho$.
The proof proceeds by constructing a canonical
homomorphism
$$\la:K_0(\LaG) \to \prod_{\rho \in \cV}K_0(L),
\eqno{(27)}$$
which will be the composition
$\la=\la_4 \circ \la_3 \circ \la_2 \circ \la_1$
of four natural maps $\la_i$ $(i=1, \dots, 4)$,
defined as follows.
Firstly, $\la_1$ is the canonical map appearing
in Lemma 3.5.
Secondly, we take
$$\la_2:K_0(\Zp[\Del])\to K_0(\Qp[\Del])$$
and
$$\la_3:K_0(\Qp[\Del])\to K_0(L[\Del])$$
to be the maps induced by the evident inclusions
of rings.
Finally, we take $\la_4$ to be the isomorphism
$$\la_4:K_0(L[\Del])\overset{\simeq}{\to}
\prod_{\rho \in \cV}K_0(M_{n_{\rho}}(L))
\overset{\simeq}{\to}
\prod_{\rho \in \cV}K_0(L),$$
where the first map is induced by $\gamma$,
and the second is given by Morita equivalence.
Now $\la_1$ is injective by Lemma 3.5.
Moreover, it is well known from the representation theory of
finite groups that $\la_2$ is injective
(see \cite{Se}, Chap. 16, Theorem 34, Corollary 2)
and that $\la_3$ is injective
(see \cite{Se}, Chap. 14, \S 14.6).
Hence we conclude that the homomorphism $\la$ is always injective.\par

To complete the proof of Proposition 3.4,
we shall employ an alternative description of the map $\la$.
Consider the map
$$\tau:K_0(\LaOG) \to K_0(O) \eqno{(28)}$$
induced by the augmentation map from $\LaOG$ to $O$.
Recall that, since $G$ has no element of order $p$, $\LaOG$ has
finite global dimension, and we can identify $K_0(\LaOG)$
with the Grothendieck group of the category of
all finitely generated $\LaOG$-modules
(see \cite{MR}, \S 12.4.8).
Let $U$ (resp. $A$) be a finitely generated $\LaOG$-module
(resp. $O$-module), and write
$[U]$ (resp. $[A]$) for the class of $U$ (resp. $A$) in $K_0(\LaOG)$
(resp. $K_0(O)$).
Then it is easily seen that
$\tau$ is given explicitly by
$$\tau([U])=\sum_{i \geq 0}(-1)^i[H_i(G, U)]. \eqno{(29)}$$
Let us also note that $\tau$ clearly factors through the map
$$\vare:K_0(\LaOG) \to K_0(\LaOGa) \eqno{(30)}$$
by the natural surjection from $\LaOG$ to
$\LaOGa$.
Moreover, $\vare$ is given explicitly by
$$\vare([U])=\sum_{i \geq 0}(-1)^i[H_i(H, U)]. \eqno{(31)}$$
Let $j:K_0(O)\to K_0(L)$ be the isomorphism induced by
the inclusion of $O$ in $L$.
For each $\rho$ in $\cV$, let $\twrho(M)$ be the $\LaOG$-module
defined by (13). For each finitely generated
$\LaG$-module $M$, it can be readily verified that
$$\la([M])=\prod_{\rho \in \cV}j(\tau([\twrho(M)])). \eqno{(32)}$$
We can now finish the proof of Proposition 3.4.
Suppose that $M$ belongs to $\MSaG$.
By Lemmas 3.1 and 3.2, we conclude that
$H_i(H, \twrho(M))$ is $\LaOGa$-torsion for all $i \geq 0$,
and thus their class in $K_0(\LaOGa)$ vanishes since the latter is
isomorphic to $\Z$, via the map which assigns to a finitely
generated $\LaOGa$-module its rank.
Hence it follows from (31) that
$\vare([\twrho(M)])=0$, whence certainly
$\tau([\twrho(M)])=0$,
for all $\rho$ in $\cV$.
But then (32) implies that $\la([M])=0$, and so we must have
$[M]=0$, because $\la$ is injective.
It now follows from the exactness of (24) that $\partial_G$
is surjective, and the proof of Proposition 3.4 is complete.
\qed\par

\bigskip

In view of Proposition 3.4,
we can now define the notion
of a {\it characteristic element} for each $M$ in
$\MSaG$.

\bigskip

{\bf {Definition.}}
For each $M$ in $\MSaG$, a characteristic
element of $M$ is any $\xi_M$ in $K_1(\LaGSa)$ such that
$$\partial_G(\xi_M)=[M]. \eqno{(33)}$$

\bigskip

We recall that we say an $\LaG$-module $M$
has finite $G$-Euler characteristic if $H_i(G,M)$
is finite for all $i \geq 0$.
If $M$ has finite $G$-Euler characteristic, we define
$$\chi(G,M)=\prod_{i \geq 0} \sharp (H_i(G,M))^{(-1)^i},
\eqno{(34)}$$
the product on the right is finite beause of our hypothesis that
$G$ has no element of order $p$.
The principal result of this section,
which is directly inspired by a parallel result in
\cite{Ve1}, is the following relation between
characteristic elements and $G$-Euler characteristics.
We write $| \ |_p$ for the valuation of
$\barQp$, normalized so that
$|p|_p=1/p$. For our continuous homomorphism
$$\rho:G \to GL_n(O),$$
we write $\mrho=[L:\Qp]$, where $L$ is the
quotient field of $O$. Let $\hatrho$ denote the contragradient
representation of $G$, i.e.
$\hatrho(g)=\rho(g^{-1})^t$ for $g$ in $G$,
where $'t\,'$ denotes the transpose matrix.

\begin{thm}
\label{t3.6}
Assume that $G$ has no element of order $p$.
Take $M$ in $\MSaG$, and let $\xi_M$ be a
characteristic element of $M$. Then,
for every continuous homomorphism $\rho:G \to GL_n(O)$
such that $\chi(G, \twhatrho(M))$ is finite,
we have
$$\xi_M(\rho)\ne0, \; \infty, \eqno{(35)}$$
and
$$\chi(G,\twhatrho(M))=|\xi_M(\rho)|_p^{-\mrho}, \eqno{(36)}$$
where $\mrho=[L:\Qp]$, and $L$ is the quotient
field of $O$.\par
\end{thm}

Before giving the proof of Theorem 3.6,
we recall an important ingredient in it,
which was first introduced in \cite{CSS}.
Assume $M$ lies in $\MSaG$.
By Lemma 3.1, the $H_i(H,M)$ $(i \geq 0)$ are finitely
generated torsion $\LaGa$-modules, which are zero for $i\geq d$.
Let $f_{i,M}$ denote a characterisitic power series for
$H_i(H,M)$ as a $\LaGa$-module.
We then define the {\it Akashi series}
$\Ak(M)$ by
$$\Ak(M)=\prod_{i \geq 0} f_{i,M}^{(-1)^i}
\ {\rm {mod}} \ \LaGa^{\times}. \eqno{(37)}$$
As is explained in \cite{CSS}, $\Ak$ induces
in the evident fashion a homomorphism
$$\Ak:K_0(\MSaG)\to Q(\Ga)^{\times}/\LaGa^{\times}.
\eqno{(38)}$$
Suppose now that $M$ is also an $O$-module,
so that we can regard the $H_i(H,M)$ as finitely
generated torsion $\LaOGa$-modules.
Let $g_{i,M}$ denote a characteristic power series
of $H_i(H,M)$ as a $\LaOGa$-module.
We then define
$$\Ak_O(M)=\prod_{i \geq 0} g_{i,M}^{(-1)^i}
\ {\rm {mod}} \ \LaOGa^{\times}. \eqno{(39)}$$
Now $\LaOGa$ is a free $\LaGa$-module of rank
$\mrho=[L:\Qp]$, and we let $N$ denote the norm
map from $\LaOGa$ to $\LaGa$.
It follows from \cite{BCA},
Chap.VII, \S 4.8, Prop. 18 that
$N(g_{i,M})=f_{i,M}$ mod $\LaGa^{\times}$ for all $i \geq 0$,
and so we conclude that
$$N(\Ak_O(M))=\Ak(M). \eqno{(40)}$$

We now begin the proof of Theorem 3.6.
Since the map $M \mapsto \twhatrho(M)$ preserves
exact sequences, we can define a map
$$\Delrho:K_0(\MSaG)\to \QOGa^{\times}/\LaOGa^{\times}\eqno{(41)}$$
by
$$\Delrho([M])=\Ak_O(\twhatrho(M)).\eqno{(42)}$$
Let
$$\partial_{\Ga}:\QOGa^{\times}\to \QOGa^{\times}/\LaOGa^{\times}$$
be the natural surjection.
Consider the diagram
$$\begin{CD}
K_1(\LaGSa)@>{\partial_G}>>K_0(\MSaG)\\
@V{\Sar}VV @VV{\Delrho}V\\
\QOGa^{\times}@>{\partial_{\Ga}}>>\QOGa^{\times}/\LaOGa^{\times}.
\end{CD}
\eqno{(43)}
$$

\begin{lem}
\label{l3.7}
The diagram {\rm (43)} is commutative.
\end{lem}

Let us now assume Lemma 3.7, and show that Theorem 3.6 follows.
Let $M$ and $\xi_M$ be as in Theorem 3.6,
and assume that $\chi(G,\twhatrho(M))$ is finite.
As before, let $\vphi:\LaOGa \to O$ be the augmentation map,
and $\fp$ its kernel. We again write
$$\vphi:\LaOGa_{\fp} \to L$$
for the homomorphism induced by $\vphi$.
Since $\chi(G,\twhatrho(M))$ is finite, it is proven
in \cite{CSS} that
$\vphi(f_{i,\twhatrho(M)})\ne 0$ for all $i \geq 0$,
that $\vphi(\Ak(\twhatrho(M))) \ne 0, \; \infty$,
and that
$$\chi(G,\twhatrho(M))=|\vphi(\Ak(\twhatrho(M)))|_p^{-1}.
\eqno{(44)}$$
Using (40)
and the obvious fact that $N$ commutes with
the augumentation map,
it follows that
$\vphi(\Ak_O(\twhatrho(M)))\ne 0, \; \infty$,
and that
$$\chi(G,\twhatrho(M))=|\vphi(\Ak_O(\twhatrho(M)))|_p^{-\mrho}.
\eqno{(45)}$$
But the commutativity of (43) shows that
$$\partial_{\Ga}(\Sar(\xi_M))=\Ak_O(\twhatrho(M)).
\eqno{(46)}$$
Hence assertions (35) and (36) are clear from (45) and (46).

\bigskip

We now prove the commutativity of the diagram (43),
which will complete the proof of Theorem 3.6.
Let $\Ktor_0(M_n(\LaOGa))$ be the Grothendieck
group of the category of all finitely
generated torsion $M_n(\LaOGa)$-modules.
Now Morita equivalence and the functoriality of
the localization sequence in $K$-theory
shows that we have a commutative diagram
$$
\begin{CD}
K_1(\LaGSa)@>{\partial_G}>>K_0(\MSaG)\\
@V{u_1}VV @VV{u_0}V\\
K_1(M_n(\QOGa))@>{\partial_n}>>
\Ktor_0(M_n(\LaOGa))\\
@V{v_1}V{\cong}V @V{\cong}V{v_0}V\\
K_1(\QOGa)@>{\partial_1}>>\Ktor_0(\LaOGa),
\end{CD}
\eqno{(47)}
$$
where $u_0$ is induced by the ring homomorphism
$\Phirho$ given in (15), $u_1$ is induced by the extension
(17) of $\Phirho$,
$\partial_n$ and $\partial_1$ are
connecting homomorphisms of localization,
and $v_0$ and $v_1$ are given by Morita equivalence.
Granted the canonical identifications
$$K_1(\QOGa)=\QOGa^{\times}, \quad
\Ktor_0(\LaOGa)=\QOGa^{\times}/\LaOGa^{\times},$$
it therefore suffices to compute $v_0\circ u_0$,
and show that it is indeed $\Delta_{\hatrho}$
given by $\hatrho(g)=\rho(g^{-1})^t$.
Take $M$ in $\MSaG$. We can find a finite resolution
$$0 \to P_r \to P_{r-1} \to \cdots \to
P_0 \to M \to 0, \eqno{(48)}$$
where each $P_i$ is a finitely generated projective
$\LaG$-module.
Thus, viewing $M_n(\LaOGa)$ as a right $\LaG$-module via
the homomorphism $\Phirho$, we have, by definition,
$$u_0([M])=\sum_{i \geq 0}(-1)^i
[M_n(\LaOGa)\otimes_{\LaG}P_i].
\eqno{(49)}$$
On the other hand, for an arbitrary ring $R$, the Morita equivalence 
between the
category of $M_n(R)$-modules and the category of $R$-modules is given by $N
\mapsto (R^n)^t \otimes_{M_n(R)} N$ where $(R^n)^t$ means the space of row
vectors over $R$ with $n$ entries ($'t\,'$ denotes the transpose) on which 
$M_n(R)$
acts from the right in the evident way.
Put $R=\LaOGa$. We endow $R\otimes_O(O^n)^t$ with
the right $G$-action given by
$(\tau \otimes x^t)g=\tau\bar{g}\otimes x^t\rho(g)^t$,
where $\tau$ is in $R$, $x$ is in $O^n$,
$g$ is in $G$, and $\bar{g}$ denotes the image of $g$ in $\Ga$.
Recalling that the left $G$-action on
$\twhatrho(P_i)=P_i\otimes_{\Zp}O^n$ is given by
$g(y \otimes x)=gy \otimes \hatrho(g)x$
for $y$ in $P_i$,
we obtain an isomorphism of $R$-modules
$$(R \otimes_O (O^n)^t)\otimes_{\LaG}P_i
\overset{\simeq}{\to}R\otimes_{\LaOG}\twhatrho(P_i)$$
by mapping $\tau\otimes x^t \otimes y$ to $\tau \otimes y \otimes x$.
Combining this isomorphism with the natural identification
$$((R^n)^t \otimes_{M_n(R)}M_n(R)) \otimes_{\LaG}P_i
=(R\otimes_O (O^n)^t) \otimes_{\LaG}P_i,$$
we deduce finally that we have an isomorphism of
$R$-modules

$$((R^n)^t \otimes_{M_n(R)}M_n(R)) \otimes_{\LaG}P_i
\overset{\simeq}{\to} R\otimes_{\LaOG}\twhatrho(P_i).$$
Hence
$$v_0([M_n(\LaOGa)\otimes_{\LaG}P_i])
=\LaOGa\otimes_{\LaOG}\twhatrho(P_i).
\eqno{(50)}$$
But it has already been remarked that,
if $N$ is any free $\LaOG$-module, then
$\twhatrho(N)$ is again free.
Hence it follows easily that
$\twhatrho(P_i)$ is again a projective $\LaOG$-module,
and the exact sequence
$$0 \to \twhatrho(P_r) \to \twhatrho(P_{r-1})
\to \cdots \to \twhatrho(P_0) \to \twhatrho(M) \to 0$$
is a projective resolution of $\twhatrho(M)$.
Tensoring this exact sequence with $\LaOGa$
over $\LaOG$,
we get a complex whose cohomology groups are
the $H_i(H, \twhatrho(M))$ $(i \geq 0)$.
Thus it follows that the image of
$v_0\circ u_0([M])$ under the canonical isomorphism
from $\Ktor_0(\LaOGa)$ to
$\QOGa^{\times}/\LaOGa^{\times}$ is precisely
$\Ak_O(\twhatrho(M))$.
This completes the proof of Theorem 3.6.
\qed\par

\bigskip

We say our continuous representation
$\rho:G\to GL_n(O)$ is an
{\it Artin representation} if $\Ker(\rho)$
is open in $G$, or equivalently if $\rho$ factors through
a finite quotient of $G$.
The following stronger form of Theorem 3.6, but for a much
more restricted class of modules,
will be needed in the last section to study some of
the consequences of the "main conjecture".

\begin{thm}
\label{t3.8}
Assume that $G$ has no element of order $p$.
Let $M$ be a module in $\MSaG$, which,
in addition, satisfies

$$H_i(H',M)=0 \quad {\text {for all  $i \geq 1$}},
\eqno{(51)}$$
and for all open subgroups $H'$ of $H$, which are normal
in $G$. Let $\xi_M$ denote a characteristic element of $M$.
Then $\xi_M(\rho)\ne \infty$ for every Artin representation
$\rho$ of $G$. Moreover, for each Artin representation
$\rho$ of $G$, $\xi_M(\rho)\ne 0$ if and only if
$\chi(G, \twhatrho(M))$ is finite.
\end{thm}

The following lemma is the essential ingredient in the proof
of Theorem 3.8.\par

\begin{lem}
\label{l3.9}
Let $M$ be a module in $\MSaG$, satisfying ${\rm {(51)}}$.
Then, if $\rho$ is any Artin representation of $G$,

$$\AkO(\twrho(M)) \in
\LaOGa[\frac{1}{p}] \;
{\rm {modulo}} \;
\LaOGa^{\times}.
\eqno{(52)}$$
Moreover, if $\vphi:\LaOGa_{\fp}\to L$ is the homomorphism
induced by the augmentation map, then
$\vphi(\AkO(\twrho(M)))\ne 0$ if and only if
$\chi(G, \twrho(M))$ is finite.
\end{lem}

\begin{pf}
Let $\rho$ be any Artin representation of $G$, and
put $W=\twrho(M)$.
Let $g_{i,W}$ denote the characteristic power series of
$H_i(H,W)$ as a $\LaOGa$-module $(i \geq 0)$. Let $\pi$ denote
a local parameter of $O$. We now prove that we can take
$$g_{i,W}=\pi^{\mu_i} \quad (i \geq 1)
\eqno{(53)}$$
for some integer $\mu_i \geq 0$. Take $H'=\Ker(\rho) \cap H$.
Hence (51) holds for $H'$, whence it follows easily
that we have $H_i(H',\MO)=0$ for all $i \geq 1$.
But since $H' \subset \Ker(\rho)$,
we have $W=\MO^n$ as $H'$-modules, and so
we conclude that
$$H_i(H',W)=0 \quad (i \geq 1).
\eqno{(54)}$$
Put $\Del=H/H'$. It follows from (53) and the Hochschild-Serre
spectral sequence that
$$H_i(H,W)=H_i(\Del,W_{H'}) \quad
(i \geq 1). \eqno{(55)}$$
But the group on the right is clearly annihilated by the order
of the finite group $\Del$, whence
it is plain that the characteristic power series of the
module on the left must be a power of $\pi$. This proves
(53). Hence
$$\AkO(W)=g_{0,W} \prod_{i \geq 1}\pi^{(-1)^i\mu_i}
\quad {\text {mod $\LaOGa^{\times}$}},
\eqno{(56)}$$
and so assertion (52) is proven.
It also follows from (56) that $\vphi(\AkO(W))\ne 0$
if and only if $\vphi(g_{0,W}) \ne 0$.
Now a standard argument with the Hochschild-Serre
spectral sequence (see \cite{CSS}) shows that
$\chi(G,M)$ is finite if and only if $\vphi(g_{i,W})\ne 0$
for all $i \geq 0$. But, in view of (53),
it is clear that $\vphi(g_{i,W}) \ne 0$ for all $i \geq 0$
if and only if $\vphi(g_{0,W}) \ne 0$. This completes
the proof of Lemma 3.9.
\end{pf}

Theorem 3.8 is an almost immediate consequence of Lemmas 3.7 and 3.9.
By definition, for any Artin representation $\rho$ of $G$,
$\xi_M(\rho) \ne \infty$ means that
$\Sar(\xi_M)$ belongs to $\LaOGa_{\fp}$.
But, by Lemma 3.7,
$$\partial_{\Ga}(\Sar(\xi_M))=\AkO(\twhatrho(M)),
\eqno{(57)}$$
whence it is clear from (52) that $\Sar(\xi_M)$ does
belong to $\LaOGa_{\fp}$.
The final assertion of Theorem 3.8 is also clear from (57)
and the final assertion of Lemma 3.9.
This completes the proof of Theorem 3.8.
\qed\par

\bigskip

We next establish an analogue of the classical Artin formalism
for $G$-Euler characteristics.
The result is of interest in its own right,
but we shall also use it to study the numerical example
at the end of this section.
In fact, as we shall  show, the Artin formalism
holds for all compact $p$-adic Lie groups $G$ with no element of
order $p$, and we do not need for this result the existence
of a closed normal subgroup $H$ of $G$ with $G/H \cong \Zp$.
We use similar notation to earlier.
Let $G'$ be an arbitrary open normal subgroup of $G$,
and let $\Del=G/G'$. Let $\cV=\cV(\Del)$
denote the set of irreducible representations
of $\Del$ over $\barQp$. We take $L$ to be any finite extension of
$\Qp$ such that all representations in $\cV$ can be realized over $L$,
and we write $O$ for the ring of integers of $L$.

\begin{thm}
\label{t3.10}
Let $G$ be any compact $p$-adic Lie group with no element of
order $p$, $G'$ an open normal subgroup of $G$,
and let $\Del=G/G'$.
Let $M$ be a finitely generated $\LaG$-module.
Then $\chi(G',M)$ is finite if and only if
$\chi(G, \twrho(M))$ is finite for all $\rho$ in
$\cV=\cV(\Del)$.
When $\chi(G',M)$ if finite, we have
$$\chi(G',M)^{[L:\Qp]}=\prod_{\rho \in \cV}
\chi(G, \twrho(M))^{n_{\rho}},
\eqno{(58)}$$
where $n_{\rho}$ is the dimension of $\rho$.
\end{thm}

We now prove Theorem 3.10.
Put $R=O[\Del]$. If $\rho$ is in $\cV$, we let $L_{\rho}$
denote a free $O$-module of rank $n_{\rho}$, endowed with a left
action of $R$ which realizes $\rho$. We then have an exact
sequence of left $R$-modules

$$0 \to R \to \bigoplus_{\rho \in \cV}L_{\rho}^{n_{\rho}}
\to W \to 0,
\eqno{(59)}$$
where $W$ is finite. If $U$ is any left $R$-module, we can view
$U$ as a right $R$-module by defining $u.\del$ to be
$\del^{-1}.u$ for $\del$ in $\Del$ and $u$ in $U$.
In particular, we can view all the modules appearing in (59)
to be right $R$-modules or $\LaOG$-modules
in this fashion.
Now take $M$ to be the module appearing in Theorem 3.10,
and put $\MO=M\otimes_{\Zp}O$.
We take a resolution

$$0 \to P_n \to P_{n-1}\to \cdots \to P_0 \to \MO
\to 0, \eqno{(60)}$$
where $P_i$ $(0 \leq i \leq n)$
is a finitely generated projective $\LaOG$-module.
Tensoring (59) on the right over $\LaOG$ with $P_i$,
we obtain an exact sequence

$$0 \to R\otimes_{\LaOG}P_i \to
\bigoplus_{\rho \in \cV}(L_{\rho}\otimes_{\LaOG}P_i)^{n_{\rho}}
\to W \otimes_{\LaOG}P_i \to 0,
\eqno{(61)}$$
the injectivity on the left is valid because
$P_i$ is projective.
Hence we obtain an exact sequence of complexes

$$0 \to C \to \bigoplus_{\rho \in \cV}D_{\rho}^{n_{\rho}}
\to K \to 0,
\eqno{(62)}$$
where
$C=(C_i)$, $D_{\rho}=(D_{\rho,i})$,
$K=(K_i)$, with $i$ running from $0$ to $n$,
are given by

$$C_i=R \otimes_{\LaOG}P_i,
\quad D_{\rho, i}=L_{\rho}\otimes_{\LaOG}P_i, \quad
K_i=W\otimes_{\LaOG}P_i.
\eqno{(63)}$$
Now
$$R\otimes_{\LaOG}\MO=(\MO)_{G'},\quad
L_{\rho}\otimes_{\LaOG}\MO\overset{\simeq}{\to}
(\twrho(M))\otimes_{\LaOG}O=(\twrho(M))_G. \eqno{(64)}$$
Writing $H_i(C)$, $H_i(D_{\rho})$, $H_i(K)$ for
the $i$-th homology groups of the complexes $C$, $D_{\rho}$, $K$,
we therefore have, for $0 \leq i \leq n$,

$$H_i(C)=H_i(G', \MO), \quad
H_i(D_{\rho})=H_i(G, \twrho(M)).
\eqno{(65)}$$
Now, as $O$ is flat over $\Zp$,

$$H_i(G',\MO)=H_i(G',M)\otimes_{\Zp}O. \eqno{(66)}$$
The long exact sequence for the cohomology of the
exact sequence of complexes (62) shows that the conclusions
of Theorem 3.10 will all follow provided we can prove that

$$
{\text {$H_i(K)$ is finite for $0 \leq i \leq n$ and }}
\prod_{i=0}^n\sharp (H_i(K))^{(-1)^i}=1.
\eqno{(67)}$$
To establish (67), we note that each $C_i$
is a projective finitely generated $R$-module,
and hence $C$ defines a class in $K_0(R)$,
which we denote by $[C]$. In fact, we claim that
$[C]=0$. Indeed, the natural map from
$K_0(R)$ to $K_0(L[\Del])$ is injective
(see \cite{Se}, Chap. 16, Theorem 34, Corollary 2),
and $[C]$ must be sent to zero under this map
because of the fact that
$H_i(C)$ is finite for $0\leq i \leq n$.
As $[C]$ maps to zero under an injective map,
we must have $[C]=0$. We now define a canonical homomorphism

$$\theta:K_0(R) \to \Qp^{\times}.
\eqno{(68)}$$
Note that, if $A$ is any finitely generated projective $R$-module,
then $W\otimes_{R}A$ is clearly finite.
Moreover, as $A$ is projective,
$\Tor_1^A(W,A)=0$, and we therefore can define
$\theta$ by sending $[A]$ to $\sharp(W\otimes_RA)$.
In particular, as $[C]=0$, we conclude that
$\theta([C])=1$.
But
$$K_i=W\otimes_RC_i \quad(i=0, \cdots, n),$$
and so $K_i$ is finite, and

$$\theta([C])=\prod_{i=0}^n\sharp(K_i)^{(-1)^i}
=\prod_{i=0}^n\sharp(H_i(K))^{(-1)^i}.
\eqno{(69)}$$
Since $\theta([C])=1$, (67) follows from (69).
This completes the proof of Theorem 3.10.
\qed\par

\bigskip

{\bf {Example.}}
We end this section by discussing a numerical example of our
algebraic theory,
arising from the arithmetic of elliptic curves.
Let $E=X_1(11)$ be the elliptic curve over $\Q$ of
conductor $11$

$$y^2+y=x^3-x^2. \eqno{(70)}$$
Take $p=5$, and let $E_{5^{\infty}}$ denote the group of
all $5$-power division points on $E$.
Define

$$\Fin=\Q(E_{5^{\infty}}). \eqno{(71)}$$
By the Weil pairing, $\Fin$ contains
$\Q(\mu_{5^{\infty}})$, where
$\mu_{5^{\infty}}$ denotes the group of all $5$-power roots
of unity. In particular, $\Fin$ contains
the cyclotomic $\Z_5$-extension of $\Q$, which we denote
by $\Qcyc$. Define the Galois groups

$$G=G(\Fin/\Q), \quad H=G(\Fin/\Qcyc), \quad
\Ga=G(\Qcyc/\Q), \eqno{(72)}$$
providing an example of compact $5$-adic Lie groups
to which our general theory applies.
We remark in passing that not only is $G$ open in $GL_2(\Z_5)$,
but in fact it is well known (\cite{Fi}, Prop. 1.5)
to be isomorphic
to the subgroup of $GL_2(\Z_5)$ consisting of all elements
$\begin{pmatrix}a&b \\c&d \end{pmatrix}$ with
$$\begin{pmatrix}a&b \\c&d \end{pmatrix}
\equiv
\begin{pmatrix}1&0 \\0&* \end{pmatrix}
\; {\rm {mod}} \; 5.$$
We shall be interested in the following two irreducible
Artin representations $\rho_i$ $(i=1,2)$ of $G$.
Let $E_2$ denote the unique elliptic curve over $\Q$ of
conductor $11$ with $E_2(\Q)=0$
(thus $E_2$ is the curve $A_2$ of \cite{Fi}).
If we write $E_5$ and $E_{2,5}$ for the groups of $5$-division points on
$E$, $E_2$, respectively,
and define
$$K_1=\Q(E_5), \quad K_2=\Q(E_{2,5}), \eqno{(73)}$$
then $K_1$ and $K_2$ are both cyclic extensions
of degree $5$ of the field $\Q(\mu_5)$, and they are
both contained in $\Fin$.
Let $\chi_i$ $(i=1,2)$ denote any non-trivial character of the
Galois group of $K_i$ over $\Q(\mu_5)$. Then we define
$\rho_i$ to be the character of the Galois group of $K_i$ over $\Q$
which is induced by $\chi_i$, and we can then view $\rho_i$ as an
Artin representation of $G$. It is easily seen that
$\rho_i$ is irreducible of degree $4$, and, in fact,
$\rho_i$ can be realized as a $4$-dimensional representation
even over $\Q$. Also $\hatrho_i=\rho_i$ for $i=1,2$.\par

We write $X(E/\Fin)$ for the compact Pontrjagin dual of the
Selmer group of $E$ over $\Fin$
(see \S 4 for a fuller discussion of Selmer groups).
Then $X(E/\Fin)$ is a finitely generated $\LaG$-module, and
it is proven in \cite{CH} that
$X(E/\Fin)$ is finitely generated over
$\LaH$. Hence $X(E/\Fin)$ belongs to our
category $\MSaG$. We write
$\xi_{E,5}$ for any choice of a characteristic element of
$X(E/\Fin)$.

\begin{prop}
\label{p3.11}
Put $M=X(E/\Fin)$. For $i=1,2$,
$\chi(G,\tw_{\rho_i}(M))$ is finite, and, in fact,
$$\chi(G,\tw_{\rho_1}(M))=5^3, \quad
\chi(G,\tw_{\rho_2}(M))=5. \eqno{(74)}$$
\end{prop}

In view of Proposition 3.11, we deduce from
Theorem 3.6 that
$\xi_{E,5}(\rho_i)\ne 0, \infty$ $(i=1,2)$ and
$$|\xi_{E,5}(\rho_1)|_5^{-1}=5^3, \quad
|\xi_{E,5}(\rho_2)|_5^{-1}=5. \eqno{(75)}$$

We now prove Proposition 3.11, by combining Theorem 3.10
with the following explicit Euler characteristic
formula for $X(E/\Fin)$ (which is Theorem 1.1 of
\cite{CH} applied to this example for the prime
$p=5$).
Put $F=\Q(\mu_5)$, and let $K$ be any finite extension of $F$,
which is contained in $\Fin$. Write

$$G_K=G(\Fin/K).$$
For each finite place $v$ of $K$, we write
$d_v$ for the
degree of $K$ completed at
$v$ over the completion of $\Q$ at $v$,
$k_v$ for the residue field, and
$\tilde{E_v}$ for the reduction of $E$ modulo $v$.
Let $\Sha(E/K)$ be the Tate-Shafarevich group of $E$ over $K$.
Assume now that

$${\text {$E(K)$ and $\Sha(E/K)(5)$ are both finite.}}
\eqno{(76)}$$
Since $X(E/\Fin)$ is finitely generated over
$\LaH$, it follows that the hypotheses of Theorem 1.1 of
\cite{CH} hold for $E$ over $K$ and the prime $p=5$.
Hence Theorem 1.1 of \cite{CH} shows that,
putting $X=X(E/\Fin)$, we have
$\chi(G_K,X)$ is finite, and is given by the
arithmetic formula

$$\chi(G_K,X)=
\frac{\sharp(\Sha(E/K)(5))}{\sharp(E(K)(5))^2}
\times \prod_{v|11}(5|d_v|_5^{-1})\times
\prod_{v|5}\sharp(\tilde{E_v}(k_v)(5))^2 \ ;
\eqno{(77)}$$
here the first product on the right of (77)
is taken over the primes $v$ of $K$ dividing 11,
and the second over primes $v$ of $K$ dividing 5.
When $K=F$, it is well known
(see \cite{CS}) that

$$E(F)(5)=\Z/5\Z, \quad
\Sha(E/F)(5)=0.$$
Since there are four primes of $F$ above 11 with $d_v=1$,
and since there is a unique prime above 5
with residue field $\F_5$, we conclude from
(77) that

$$\chi(G_F,X)=5^4.\eqno{(78)}$$
Taking next $K=K_1$, it is proven in \cite{Fi} that

$$E(K_1)=(\Z/5\Z)^2, \quad
\Sha(E/K_1)(5)=(\Z/5\Z)^2,$$
and that there are four primes of $K_1$ above 11 with $d_v=5$,
and five primes of $K_1$ above 5 with residue field $\F_5$.
Hence we conclude from (77) that

$$\chi(G_{K_1}, X)=5^{16}.
\eqno{(79)}$$
Applying Theorem 3.10 to both of the open subgroups
$G'=G_F$ and $G'=G_{K_1}$, we deduce that
$\chi(G,\tw_{\rho_1}(X))$ is finite, and

$$\chi(G, \tw_{\rho_1}(X))^4=5^{12}, $$
which proves the first assertion of Proposition 3.11.
Now take $K=K_2$. It is proven in \cite{Fi} that

$$E(K_2)=\Z/5\Z, \quad
\Sha(E/K_2)(5)=0,$$
and that there are four primes of $K_2$ above 11
with $d_v=5$, and one prime of $K_2$ above 5
with residue field $\F_5$. Hence we conclude from
(77) that

$$\chi(G_{K_2},X)=5^8. \eqno{(80)}$$
Applying Theorem 3.10 to both of the open subgroups
$G'=G_F$ and $G'=G_{K_2}$, we deduce that
$\chi(G,\tw_{\rho_2}(X))$ is finite, and

$$\chi(G,\tw_{\rho_2}(X))^4=5^4,$$
which proves the second assertion of Proposition 3.11.
\qed\par

\section{Additional properties of characteristic elements.}

We now establish some additional properties of
the characteristic elements of modules in the
category $\MSaG$, and of the group
$K_1(\LaGSa)$ in which they lie.
We end this section with some conjectures about the integrality
properties of characteristic elements.\par

We begin with some general remarks.
Let $R$ be a ring with unit element,
which is left and right Noetherian.
We recall that the Jacobson radical of $R$, which we denote
by $\Jac(R)$, is the intersection of all
maximal left ideals of $R$, or equivalently the intersection
of all maximal right ideals of $R$.
We say that $R$ is {\it {semi-local}} if the ring
$R/\Jac(R)$ is both left and right Artinian.
The following lemma is well known
(see \cite{Ba}, Chap.III, Prop. 2.12,
and Chap.IX, Prop. 1.3).

\begin{lem}
\label{l4.1}
Let $I$ be a two sided ideal of a ring $R$,
and assume that $R$ is $I$-adically complete.
If $I$ is contained in $\Jac(R)$, then the natural map
from $K_0(R)$ to $K_0(R/I)$ is an isomorphism.
\end{lem}

If $G$ is an arbitrary compact $p$-adic Lie group,
it is well known (see \cite{NSW}, Chap.V, Prop. 5.2.16)
that $\LaG$ is a semi-local ring.
We assume for the rest of this section that again $G$ has a closed
normal subgroup $H$ such that $G/H$ is isomorphic to $\Zp$,
and we let $S$ be the Ore set in $\LaG$ defined in \S 2.

\begin{prop}
\label{p4.2}
The ring $\LaGS$ is semi-local.
\end{prop}

We first establish the following lemma.

\begin{lem}
\label{l4.3}
Let $J$ be any pro-$p$ open subgroup of
$H$ which is normal in $G$, and let
$\psiJ$ be the homomorphism {\rm (3)}.
Define $S_J=\psiJ(S)$. Then $S_J$ is an Ore
set of non-zero divisors in
$\OmGJ$, and the ring $\OmGJSJ$ is Artinian.
\end{lem}

\begin{pf}
Lemmas 2.1 and 2.2 show that
$S_J$ consists of all non-zero divisors in $\OmGJ$.
Moreover, as $S$ is an Ore set and $\psiJ$ is
surjective, it follows that $S_J$ satisfies the Ore condition.
Thus $\OmGJSJ$ is the total ring of
quotients of $\OmGJ$. Now choose a subgroup
$\Pi$ of $G/J$ satisfying (2), and let
$\OmPi=\Fp[[T]]$ and its quotient field
$R(\Pi)$ be as in \S 2.
Define $\Sig$ to be the set of non-zero elements of
$\OmPi$. If $\alpha$ is any element of $\Sig$,
it is clear that multiplication by $\alpha$ on
the right or the left induces an automorphism
of the $R(\Pi)$-vector space $V(G/J)$ defined by (6).
Hence, by the proof of Lemmas 2.1 and 2.2,
$\alpha$ belongs to $S_J$; in particular,
$\alpha$ is not a divisor of zero in $\OmGJ$.
As $\OmPi$ is in the centre of $\OmGJ$,
it follows that $\Sig$ is an Ore set in
$\OmGJ$.
We claim that

$$\OmGJSJ=\OmGJ_{\Sig}.  \eqno{(81)}$$
Note that (81) proves the assertion of
Lemma 4.3, because the fact that $\OmGJ$ is a
finitely generated $\OmPi$-module implies
that $\OmGJ_{\Sig}$ is a finite dimensional vector
space over $R(\Pi)$, and so $\OmGJ_{\Sig}$ is Artinian.
As $\Sig \subset S_J$, to prove
(81) we must show that every element of $S_J$ is
invertible in $\OmGJ_{\Sig}$.
Take any $s$ in $S_J$. As $\OmPi$ lies in the centre of
$\OmGJ$, which is finitely generated over $\OmPi$,
$s$ satisfies an equation (see \cite{MR}, Lemma 5.3.2)
of the form

$$s^n+a_1s^{n-1}+\cdots +a_n=0
\quad (a_i \in \OmPi). \eqno{(82)}$$
Moreover, since $s$ is not a zero divisor in
$\OmGJ$, we can assume that $a_n \ne 0$.
Now this equation can be rewritten as

$$(b_0s^{n-1}+b_1s^{n-1}+ \cdots + b_{n-1})s=1,\; $$
where $b_0=-{1/a_n}$, and
$b_i=-a_i/a_n$ $(i=1, \cdots , n-1)$,
proving that $s$ has an inverse in
$\OmGJ_{\Sig}$, as required.
This completes the proof of Lemma 4.3.
\end{pf}

We now can prove Proposition 4.2.
If $M$ is any left $\LaG$-module,
we define

$$M_S=\LaGS\otimes_{\LaG}M.$$
Let $J$ be any subgroup of $G$
as in Lemma 4.3, and let $I$ denote the kernel of
$\psiJ$. Since $\LaGS$ is flat over $\LaG$, we have
the exact sequence

$$0 \to I_S \to \LaGS \to \OmGJ_S \to 0.$$
But

$$\OmGJ_S=\OmGJSJ.$$
Hence, by lemma 4.3, it suffices to show that $I_S$
is contained in the Jacobson radical of
$\LaGS$.
In view of \cite{MR}, Theorem 0.3.8, we must prove that
$1-x$ is a unit in $\LaGS$ for every $x$ in $I_S$.
Write $x=s^{-1}\theta$ with $\theta$ in $I$ and $s$ in $S$.
Since $\theta$ is in $I$,
$\psiJ(s-\theta)=\psiJ(s)$.
Hence, by Lemma 2.1, $s-\theta$ belongs to $S$,
But then

$$1-x=s^{-1}(s-\theta)$$
is clearly a unit in $\LaGS$, and the proof of
Proposition 4.2 is complete.
\qed\par

\bigskip

We remark that $\LaGSa$ is not, in general, a
semi-local ring, as is shown by the following
example.
Take $G=H \times K$, where both $H$
and $K$ are isomorphic to $\Zp$.
Then we can identify $\LaG$ with
$R=\Zp[[U,V]]$ by mapping fixed topological generators
of $H$ and $K$ to
$U+1$ and $V+1$, respectively.
Then the set $S$ of \S 2 in this case is the complement
of the prime ideal $(p,U)$ in $R$.
Thus $R_S$ is a local ring of dimension 2, whose maximal ideal
is generated by $p$ and $U$. Hence
$\RSa=R_S[1/p]$ is a ring of dimension 1.
But $\RSa$ has infinitely many prime ideals, and is certainly
not semi-local.
Indeed, if $g$ is any
irreducible distinguished polynomial in $\Zp[U]$,
then $gR$ is a prime ideal whence $g\RSa$ is also
prime because $S^*\cap gR$ is empty.\par

\begin{thm}
\label{t4.4}
Assume that $G$ has no element of order $p$.
Then the natural maps
$\LaGS^{\times}\to K_1(\LaGS)$
and $\LaGSa^{\times} \to K_1(\LaGSa)$ are
both surjective.
\end{thm}

If $R$ is any semi-local ring,
it is well known that the natural map
$R^{\times} \to K_1R$ is surjective
(see \cite{Va}, \cite{Va2}).
In view of Proposition 4.2, this establishes the
first assertion of Theorem 4.4, and we
now give the proof of the second assertion.
We establish two preliminary lemmas.
Note that, if $P$ is a direct summand of
$\LaG$ viewed as a
left $\LaG$-module, then $P$ is a
finitely generated projective $\LaG$-module,
and $P/pP$ is a finitely generated projective
$\OmG$-module.

\begin{lem}
\label{l4.5}
The group $K_0(\OmG)$ is generated by the classes
$[P/pP]$, where $P$ runs over all direct summands of
$\LaG$ viewed as a left $\LaG$-module.
\end{lem}

\begin{pf}
Let $I=p\LaG$. Then $I$ is contained in the Jacobson radical
of $\LaG$ since,
for each $\la$ in $I$, $1-\la$ is invertible with inverse
$\sum_{n=0}^{\infty}\la^n$.
Moreover, $\LaG$ is $I$-adically complete,
and so, by Lemma 4.1, the canonical map
from $K_0(\LaG)$ to $K_0(\OmG)$ is an isomorphism.
In particular, $K_0(\OmG)$ is generated by the classes
$[P/pP]$, where $P$ ranges over all finitely
generated projective $\LaG$-modules.
To complete the proof of lemma, we must show that
we still get a set of generators if
we allow $P$ to only run over the direct
summands of $\LaG$.
Pick a pro-$p$ open normal subgroup $U$ of $G$, and
put $\Del=G/U$.
Let $D$ denote the kernel of the natural map
from $\LaG$ to $\Fp[\Del]$, so that $D$ is generated
as a $\LaG$-module by the set $W_U$ consisting of $p$ and
$1-u$, where $u$ ranges over $U$. Now $\LaU$ is a
local ring,
whose maximal ideal $\Jac(\LaU)$ is clearly generated
by $W_U$.
But $\LaG$ is finitely generated over $\LaU$, and so
it is well known that
(\cite{MR}, Prop. 9.1.3)

$$\Jac(\LaG) \cap \LaU \subseteq \Jac(\LaU).$$
This proves that $D$ is contained in $\Jac(\LaG)$,
and, as $\LaG$ is plainly $D$-adically complete,
we conclude from Lemma 4.1 that the natural map
from $K_0(\LaG)$ to $K_0(\Fp[\Del])$ is an isomorphism.
Put

$$R=\Fp[\Del]/\Jac(\Fp[\Del]).$$
Since $\Fp[\Del]$ is Artinian,
$\Jac(\Fp[\Del])$ is nilpotent,
and so Lemma 4.1 proves that the
canonical map from $K_0(\Fp[\Del])$ to $K_0(R)$
is an isomorphism. However, $R$ is a semisimple
ring, and thus $K_0(R)$ is generated by the
simple projective $R$-modules, and these are given
precisely by the direct summands of $R$ viewed as
a left $R$-module. It is then well known
(see \cite{Ba}, Chap.III, Prop. 2.12)
that one can lift these generators of
$K_0(R)$ successively back to $K_0(\Fp[\Del])$ and
to $K_0(\LaG)$,
so that they remain direct summands of
$\Fp[\Del]$ and $\LaG$, respectively. This completes
the proof of Lemma 4.5.
\end{pf}

\begin{lem}
\label{l4.6}
Put $R=\LaG[1/p]$, and assume that
$G$ has no element of order $p$.
Then the natural map $R^{\times}\to K_1R$
is surjective.
\end{lem}

\begin{pf}
Let $\cC$ denote the category of all
finitely generated $\LaG$-modules,
which are annihilated by some power of $p$.
Then, by d\'evissage, it is
well known (see, for example, \cite{MR},
Theorem 12.4.7) that we can
naturally identify $K_0(\cC)$ with $K_0(\cD)$,
where $\cD$ is the category of all finitely generated
$\OmG$-modules.
As $G$ has no element of order $p$,
$\OmG$ has finite global dimension, an so we can
identify $K_0(\cD)$ with
$K_0(\OmG)$.
Thus the localization sequence of $K$-theory for the multiplicative
Ore set consisting of the powers of $p$ lying in $\Z$
gives the commutative diagram

$$\begin{CD}
\LaG^{\times}
@>{\hookrightarrow}>>
R^{\times}
@>g>>K_0(\OmG)\\
@VVV @VfVV \Vert\\
K_1(\LaG) @>>> K_1(R) @>{\partial}>> K_0(\OmG)
\end{CD}
\eqno{(83)}$$
where $R=\LaG[1/p]$,
$f$ is the natural map,
and the bottom row is exact.
Also $\partial$  is the connecting homomorphism,
and $g=\partial \circ f$.
Since $\LaG$ is semi-local, the left vertical
map is surjective (see \cite{Va}, \cite{Va2}).
We conclude from (83) that $f$ will be surjective provided
we can prove $g$ is surjective.
We now proceed to show this.\par

By Lemma 4.5, $K_0(\OmG)$ is generated by the classes
$[P/pP]$, where $P$ runs over the direct summands of
$\LaG$, viewed as a left $\LaG$-module. Hence it suffices
to show that each such class belongs to the image of $g$.
Take such a $P$, and let $P'$ be the
$\LaG$-module such that

$$P \oplus P'=\LaG. \eqno{(84)}$$
In terms of the decomposition (84), we can define
endomorphisms of $\LaG$, viewed as a left
$\LaG$-module,
by

$$\alpha(x+y)=px+y, \quad
\beta(x+y)=x+py
\qquad (x \in P, \; y \in P').$$
Thus $\alpha$ and $\beta$ are given,  respectively,
by multiplication on the right by $a=\alpha(1)$,
$b=\beta(1)$.
Since $ab=p$, it follows
that $a$ and $b$ are both units in $R$.
On the other hand, the explicit description of
$\partial$ in the localization sequence (see \cite{Sw})
shows that

$$\partial(f(a))=[\Coker \alpha]=[P/pP].$$
This completes the proof of Lemma 4.6.
\end{pf}

We can now prove the second assertion of
Theorem 4.4.
As before, we let $R=\LaG[1/p]$,
$\OmG=\LaG/p\LaG$,
and put

$$R'=\LaGS[\frac{1}{p}]=\LaGSa, \quad
\OmGp=\LaGS/p\LaGS.$$
As in the proof of Lemma 4.6, we can identify
the Grothendieck group of the category of all finitely
generated $\LaGS$-modules, which are annihilated by
some power of $p$, with $K_0(\OmGp)$.
Thus, parallel to (83), we have the
exact sequence of localization

$$K_1(\LaGS) \to K_1(R')
\overset{\partial '}{\to}
K_0(\OmGp).
\eqno{(85)}$$
Let $f'$ denote the natural map from
$(R')^{\times}$ to $K_1(R')$, and put
$g'=\partial ' \circ f'$.
As $\LaGS$ is semi-local, the
natural map from $\LaGS^{\times}$ to $K_1(\LaGS)$
is surjective.
In view of this last remark, we conclude from the
exact sequence (85) that $f'$ is surjective provided
$g'$ is surjective. To establish the surjectivity of $g'$,
we note that we have a commutative diagram

$$
\begin{CD}
R^{\times} @>g>> K_0(\OmG)\\
@VVV @VVhV\\
(R')^{\times} @>{g'}>> K_0(\OmGp),
\end{CD}
$$
where the vertical maps are the natural ones.
We claim that both $g$ and $h$ are surjective,
which will clearly imply that $g'$ is surjective.
But the surjectivity of $g$ was proven in the course of the
proof of Lemma 4.6.
Let $S'$ denote the image of $S$ in $\OmG$.
By (iii) of Lemmas 2.1 and 2.2, $S'$ consists of non-zero
divisors in $\OmG$, and it is plain that $S'$ is an Ore set
in $\OmG$.
As localization with respect to $S$ is an exact functor,
we see easily that

$$\OmGp=\OmGSp.$$
But then it is well known (\cite{MR},
Theorem 12.4.9) that the natural map $h$ from
$K_0(\OmG)$ to $K_0(\OmGSp)$ is surjective.
This completes the proof of Theorem 4.4.
\qed\par

\bigskip

We next show that any
module $M$ in the category $\MSaG$
has many twists $\rho$ by continuous representations
of $G$ of the form (12) with finite $G$-Euler characteristics.
Indeed, consider
continuous representations

$$\eta:\Ga \to O ^{\times}. \eqno{(86)}$$
Since $\Ga=G/H$, we can always view
such $\eta$ as a continuous representation of $G$.

\begin{lem}
\label{l4.7}
Assume that $G$ has no element of order $p$, and
let $M \in \MSaG$. Then,
for all but finitely many Artin characters $\eta$ of
$\Ga$, $\chi(G, \tweta(M))$ is finite.
Moreover, if $\eta$ is a fixed character of $\Ga$
of infinite order,
$\chi(G, \twetan(M))$
is finite for all but finitely many $n \in \Z$.
\end{lem}

\begin{pf}
We pick a topological generator $\ga_0$ of $\Ga$,
and identify $\LaOGa$ with
$O[[T]]$ by mapping $\ga_0$ to $1+T$.
As earlier, let $\MO=M\otimes_{\Zp}O$.
Since $M$ belongs to $\MSaG$, $H_i(H, \MO)$ is a
torsion $\LaOGa$-module for all $i \geq 0$, and
we write $f_{i,M}$ for a characteristic
series of $H_i(H,\MO)$.
Since $\eta(H)=1$, we have

$$H_i(H,\tweta(M))=H_i(H, \MO)\otimes_O A_{\eta},$$
where $A_{\eta}$ denotes a free $O$ -module of rank $1$ on
which $\Ga$ acts via $\eta$, i.e. $\ga . a=\eta(\ga)a$
for $\ga$ in $\Ga$ and $a$ in $A_{\eta}$.
It follows easily that

$$f_{i, \tweta(M)}(T)=f_{i,M}(\eta(\ga_0)^{-1}(1+T)-1).
\eqno{(87)}$$
Now a standard argument with the Hochschild-Serre
spectral sequence (see, for example, \cite{CSS})
shows that $\chi(G, \tweta(M))$ is finite provided
$f_{i, \tweta(M)}(0)\ne 0$, or equivalently

$$f_{i,M}(\eta(\ga_0)^{-1}-1)\ne 0 \qquad
(0 \leq i \leq d-1) \ ;
\eqno{(88)}$$
here $d$ denotes the dimension of $G$ as a $p$-adic
Lie group.
Let $\fZ$ denote the set of all zeros of the
$f_{i,M}(T)$ $0 \leq i \leq d-1)$ lying in the
maximal ideal of the ring of integers of
$\barQp$.
By the Weierstrass Preparation Theorem,
$\fZ$ is a finite set. Hence,
as $\eta$ runs over all the Artin characters of $\Ga$,
only finitely many of the $\eta(\ga_0)^{-1}-1$
can lie in $\fZ$, and the first assertion
of the lemma is proven.
Now suppose that $\eta$ has infinite order,
so that the characters $\eta^n$ $(n \in \Z)$ are all distinct.
Thus, for fixed $\eta$ of infinite order, only finitely
many of the $\eta(\ga_0)^{-n}-1$,
as $n$ runs over $\Z$, can lie in $\fZ$, and
the second assertion of Lemma 4.7 follows.
\end{pf}

The following two examples illustrate
what differences occur between the commutative
and the non-commutative theory,
especially when one considers modules which lie outside
the category $\MSaG$ in the non-commutative case.\par

\bigskip

{\bf {Example.}}
Suppose now that $G$ is commutative and has no element of
order $p$, and let $M$ be any finitely generated
torsion $\LaG$-module.
A slight generalization of a lemma of Greenberg \cite{Gr}
then shows that we can always find a closed subgroup $H$ of $G$
such that $\Ga=G/H$ is isomorphic to $\Zp$,
and $M$ belongs to the category $\MSaG$. Thus, in this case,
one can always find, by Lemma 4.7, continuous
representations $\rho$ of $G$ such that
$\chi(G, \twrho(M))$ is finite. The following
example shows that this last assertion is false, in general,
when $G$ is not commutative.
Take $G=GL_3(\Zp)$, and assume that
$p \geq 5$ to ensure that $G$ has no element of order $p$.
Let

$$
u=
\begin{pmatrix}
1&0&1\\
0&1&0\\
0&0&1
\end{pmatrix}, \quad
v=\begin{pmatrix}
1&0&0\\
0&1&1\\
0&0&1
\end{pmatrix}, \quad
f=u-v.
$$
Consider the module

$$M=\LaG/\LaG f.$$
We claim that, for every continuous representation
$\rho$ of $G$ of the form (12),
the Euler characteristic
$\chi(G, \twrho(M))$ is never finite.
Indeed, a free resolution of $\twrho(M)$ by
$\LaOG$-modules is given by

$$0 \to \twrho(\LaG) \to \twrho(\LaG)
\to \twrho(M) \to 0,$$
since, as was remarked earlier
(see \cite{Ve1}), $\twrho(\LaG)$
is a free $\LaOG$-module of rank $n$.
Taking $G$-homology of this complex,
it follows that
$H_i(G, \twrho(M))=0$ for
$i\geq 2$, and that

$$H_1(G, \twrho(M))=\Ker(\rho(f)),
\quad
H_0(G, \twrho(M))=\Coker(\rho(f)).$$
But we now prove that $0$ is an eigenvalue
of $\rho(f)$ for every continuous $\rho$.
We first note that for $a\in \Z_p$, we have $a$-th powers
$u^a$ and
$v^a$ defined as

$$
u^a=
\begin{pmatrix}
1&0&a\\
0&1&0\\
0&0&1
\end{pmatrix}, \quad
v^a=\begin{pmatrix}
1&0&0\\
0&1&a\\
0&0&1
\end{pmatrix}.
$$
  Taking $a$ in $\Zp^{\times}$,
we define

$$\sig_a=
\begin{pmatrix}
a&0&0\\
0&1&0\\
0&0&1
\end{pmatrix},
\quad
\tau_a=
\begin{pmatrix}
1&0&0\\
0&a&0\\
0&0&1
\end{pmatrix},
$$ and one verifies immediately that
$$\sig_au\sig_a^{-1}=u^a, \quad
\tau_a v \tau_a^{-1}=v^a.$$
Hence the eigenvalues of $\rho(u)$ and $\rho(v)$
in $\barQp$ must be stable under exponentiation by
any $a$ in $\Zp^{\times}$.
Since there are only finitely many of these eigenvalues,
they must therefore all be $p$-power roots of unity.
Moreover, $\rho(u)$ and $\rho(v)$ commute because
$uv=vu$.
Hence $\rho(u)$ and $\rho(v)$ must have a common
eigenvector $z$, with respective eigenvalues
$\alpha$ and $\beta$. Since $\alpha$ and $\beta$ are
both $p$-power roots of unity, we have either

$${\text {$\beta=\alpha^w$  or  $\alpha=\beta^w$}}
\eqno{(89)}$$
for some element $w$ in $\Z$.
Define

$$g=
\begin{pmatrix}
1&w-1&0\\
0&1&0\\
0&0&1
\end{pmatrix},
\quad
h=
\begin{pmatrix}
1&0&0\\
w-1&1&0\\
0&0&1
\end{pmatrix}.$$
One verifies that

$$g^{-1}ug=u,
\quad
g^{-1}vg=u^{1-w}v,
\quad
h^{-1}uh=uv^{1-w},
\quad
h^{-1}vh=v.
\eqno{(90)}$$
If the first option of (89) holds,
one uses the first two equations of (90) to verify that
$\rho(g)z$ is an eigenvector of both
$\rho(u)$ and $\rho(v)$ with the same eigenvalue $\alpha$.
Similarly, if the second option of (89) is true,
one uses the second two equations of (90)
to conclude that $\rho(h)z$ is an eigenvector for
both $\rho(u)$ and $\rho(v)$ with the same eigenvalue.
In either case, we see that $0$ is an eigenvalue of
$\rho(f)$, completing the proof that
$\chi(G,\twrho(M))$ is never finite.\par

\bigskip

{\bf {Example.}}
We give a second example to illustrate the differences between
the commutative and non-commutative theory.
Assume that $G$ is $p$-valued in the sense of Lazard \cite{La}
(for example, provided $p > n+1$, every pro-$p$ closed
subgroup of $GL_n(\Zp)$ is $p$-valued).
Then it follows from \cite{La} that both
$\LaG$ and $\OmG=\LaG/p\LaG$ have no zero divisors.
We consider any $\LaG$-module of the form
$M=\LaG/\LaG w$, where $w$ is an element of $\LaG$, which is
not a unit, and which does not belong to $p\LaG$.
Since $\OmG$ has no zero divisors, we see that
$M$ has no $p$-torsion.
Moreover, it is easily seen that $M$ is not pseudo-null
as a $\LaG$-module in the sense of \cite{Ve1}.
Further, in a similar manner to that of the
previous example, we see that for every continuous
representation $\rho$ of $G$,
we have $H_i(G, \twrho(M))=0$ for $i \geq 2$,
and

$$H_1(G, \twrho(M))=\Ker(\rho(w)), \quad
H_0(G, \twrho(M))=\Coker(\rho(w)).$$
Suppose first that $G$ is commutative.
Thus $G$ is isomorphic to $\Zp^d$ for some integer
$d \geq 1$, and $\LaG$ is isomorphic to the ring
$\Zp[[T_1, \cdots, T_d]]$ of formal power series
in $d$ variables with coefficients in $\Zp$.
Identifying $w$ with a formal power series $w(T_1, \cdots, T_d)$,
and recalling that $w$ is not a unit, and is not divisible
by $p$, a well known argument with the
Weierstrass preparation theorem shows that there always
exist $\alpha_1, \cdots, \alpha_d$
in the maximal ideal of the ring of integers of
$\barQp$ such that
$w(\alpha_1, \cdots, \alpha_d)=0$.
Let $O$ be the ring of integers
of any finite extension of $\Qp$
containing
$\alpha_1, \cdots, \alpha_d$. We can define
$\rho:G \to O^{\times}$
by specifying $\rho(\ga_i)=\alpha_i+1$
$(1 \leq i \leq d)$,
where $\ga_1, \cdots, \ga_d$ denote the $\Zp$-basis of
$G$ with $\ga_i=T_i+1$ $(1 \leq i \leq d)$.
Clearly $\rho(w)=0$,
and thus $\twrho(M)$ does not have finite
$G$-Euler characteristic.
On the other hand, take $G$ to be kernel of the
reduction map from $GL_3(\Zp)$ to
$\GL_3(\F_p)$, and assume that $p \geq 5$ to
ensure that $G$ is $p$-valued. Define

$$x=
\begin{pmatrix}
1&0&p\\
0&1&0\\
0&0&1
\end{pmatrix},
\quad
w=x-1-p.$$
An entirely analogous argument to that given for $u$  and $v$
in the previous example shows that, for every continuous
representation $\rho$ of $G$, the eigenvalues of $\rho(x)$
in $\barQp$ are all $p$-power roots of
unity.
Thus the eigenvalues of $\rho(w)$ in $\barQp$
are never zero, and so it follows that
$\chi(G, \twrho(M))$ is finite for every continuous
representation $\rho$ of $G$.\par

\bigskip

We end this section by making some conjectures about
the integrality properties of characteristic elements of modules
in the category $\MSaG$.
For brevity,  let us write

$$R_1=\LaG, \quad
R_2=\LaG[\frac{1}{p}], \quad
R_3=\LaGSa.
\eqno{(91)}$$
Since $S^*$ consists of non-zero divisors in
$\LaG$, we have natural inclusions
$R_1\subset R_2 \subset R_3$. We also write

$$A_1=\LaOGa, \quad A_2=\LaOGa[\frac{1}{p}].
\eqno{(92)}$$
It is also convenient to give a name $\alpha$ to the
canonical map

$$\alpha:R_3^{\times} \to K_1(R_3),
\eqno{(93)}$$
which we know is surjective by Theorem 4.4.
Finally, we write $\cRG$ for the set of all continuous
representations $\rho$ of $G$ of the form (11),
where $O$ is allowed to run over the rings of integers
of all finite extensions of $\Qp$, and
$\cAG$ for the subset of $\cRG$ consisting of the
Artin representations.
In the next conjecture, the statements for all $\rho$
in $\cRG$, and for all $\rho$ in
$\cAG$ will just be abbreviated to the respective
assertions for all $\rho$, and for all Artin $\rho$.\par

\begin{conj}
Assume that $G$ has no element of order $p$.
In each of the following four cases,
the assertions
${\rm (a)}$, ${\rm (b)}$, ${\rm (c)}$ and
${\rm (d)}$ are equivalent for
any element $\xi$ of $K_1(R_3)$:-\par

\bigskip

Case $1$.
${\rm (a)}$ $\xi \in \alpha(R_2 \cap R_3^{\times});$
${\rm (b)}$ $\xi(\rho) \ne \infty$ for all
$\rho;$
${\rm (c)}$ $\Sar(\xi) \in A_2$ for all $\rho;$
${\rm (d)}$ $\Sar(\xi) \in A_2$ for all Artin $\rho$.\par

\bigskip

Case $2$.
${\rm (a)}$ $\xi \in \alpha(R_1 \cap R_3^{\times});$
${\rm (b)}$ $\xi(\rho)$ is finite and in $O$ for all $\rho;$
${\rm (c)}$ $\Sar(\xi) \in A_1$ for all $\rho;$
${\rm (d)}$ $\Sar(\xi) \in A_1$
for all Artin $\rho$.\par

\bigskip

Case $3$.
${\rm (a)}$ $\xi \in \alpha(R_2^{\times});$
${\rm (b)}$ $\xi(\rho) \ne 0, \infty$ for all $\rho;$
${\rm (c)}$ $\Sar(\xi) \in A_2^{\times}$ for all $\rho;$
${\rm (d)}$ $\Sar(\xi) \in A_2^{\times}$ for all Artin $\rho$.\par

\bigskip

Case $4$.
${\rm (a)}$ $\xi \in \alpha(R_1^{\times});$
${\rm (b)}$ $\xi(\rho)$ is finite and in $O^{\times}$
for all $\rho;$
${\rm (c)}$ $\Sar(\xi) \in A_1^{\times}$ for all $\rho;$
${\rm (d)}$ $\Sar(\xi) \in A_1^{\times}$ for all Artin $\rho$.\par
\end{conj}

We remark that Case $1$ of Conjecture $4.8$ has an interesting
consequence for the modules $M$ in $\MSaG$
which, in addition, satisfy condition ${\rm (51)}$.
Let $M$ be any such module, and let $\xiM$ be a characteristic
element of $M$.
It follows from ${\rm (52)}$ and ${\rm (57)}$ that
$\Sar(\xiM)$ belongs to $A_2$ for all $\rho$ in $\cAG$.
Hence, if Conjecture $4.8$ is valid,
the assertions ${\rm (a)}$, ${\rm (b)}$, and ${\rm (c)}$ of
Case $1$ would be true for $\xiM$, and, in particular,
we would have $\xiM \in \alpha(R_2 \cap R_3^{\times})$.
But we have to confess that at present we are unable to prove
${\rm (a)}$, ${\rm (b)}$, and ${\rm (c)}$
of Case $1$ for the characteristic elements of such modules $M$.

\begin{lem}
In each of the four cases of Conjecture $4.8$,
${\rm (a)}$ implies ${\rm (b)}$, ${\rm (b)}$ and ${\rm (c)}$
are equivalent, and
${\rm (c)}$ implies
${\rm (d)}$.
\end{lem}

\begin{pf}
Plainly (c) implies (d), and it is also clear
that (c) implies (b) since, by definition,
$\xi(\rho)$ is the image of
$\Sar(\xi)$ under the augmentation map once
$\Sar(\xi)$ belongs to the localization of $A_1$ at the kernel
of the augmentation map.
The implication (a) implies (b) holds because if $g$ is any element of
$R_2^{\times}$ of the form $g=p^{-m}f$ with $f$ in $R_1$
and $m$ an integer $\geq 0$, then, putting
$\xi=\alpha(g)$,
$\xi(\rho)=p^{-mn}\det \rho(f)$ for every $\rho$ in $\cRG$. To prove
that (b) implies (c), we make use of the following general remark.
Recall that we have identified $A_1$ with $O[[T]]$ by
mapping a fixed topological generator $\ga_0$ of $\Ga$ to
$1+T$. Let $\eta:\Ga \to O^{\times}$ be any continuous
homomorphism. Then we claim that,
for all $\xi$ in $K_1(R_3)$ and all $\rho$ in $\cRG$,
we have the identity

$$\Phirhoetap(\xi)=\Sar(\xi)(\eta(\ga)^{-1}(1+T)-1),
\eqno{(94)}$$
where both sides of this equation are viewed as elements
of the quotient field of $O[[T]]$.
To establish (94), we note that $K_1(R_3)$ is generated as an abelian
group by those $\xi$ which are characteristic elements of modules
$M$ in $\MSaG$, and, for such $\xi$, (94) follows from (87) and
Lemma 3.7. We now use (94) to show that (b) implies (c) in Case 1.
Suppose that, for some $\rho$ in $\cRG$,
we have $\Sar(\xi)=f(T)/g(T)$, where $f(T)$ and $g(T)$ are relatively
prime elements of $O[[T]]$, and $g(T)$ is not a power of
the local parameter $\pi$ of $O$. Thus $g(T)$ must have a zero
$T=\alpha$ in the maximal ideal of the ring of integers of $\barQp$.
Enlarging $O$ if necessary, we can assume that
$\alpha$ belongs to $O$, and define a continuous homomorphism
$\eta:\Ga \to O^{\times}$ by $\eta(\ga_0)=(1+\alpha)^{-1}$. It
is then clear from (94) that $\xi(\eta \rho)=\infty$, contradicting
our assumption that (b) is valid.
Similar arguments prove that (b) implies (c) in the remaining
three cases of Conjecture 4.8,
and the proof of Lemma 4.9 is complete.
\end{pf}

The strongest evidence we have for Conjecture 4.8
is that it is true when $G$ is of the form
$\Zp^d \times \Del$, where $d$ is any integer $\geq 1$
and $\Del$ is a finite group of order prime to $p$.
However, we omit the detailed proof in this case.

\section{The main conjecture.}

For brevity, we only discuss here the main
conjecture for elliptic curves over $\Q$ over the
field generated by the coordinates of
all $p$-power division points on the curve,
where $p$ is a prime $\geq 5$ of good
ordinary reduction. However,
it will be clear that a similar conjecture can
be formulated in great generality for motives
over $p$-adic Lie extensions of number fields
which contain the cyclotomic $\Zp$-extension
of the base field.

If $F$ is a finite extension of $\Q$,
we write $\Fcyc$ for the cyclotomic $\Zp$-extension
of $F$, and put $\GaF=G(\Fcyc/F)$.
Let $E$ be an elliptic curve defined over $\Q$, and
$\Epinf$ the group of all $p$-power division points on
$E$. We define

$$\Fin=\Q(\Epinf).
\eqno{(95)}$$
By the Weil pairing, $\Fin \supset \Q(\mu_{p^{\infty}})$,
where $\mu_{p^{\infty}}$ denotes the group of all $p$-power
roots of unity. Hence $\Fin$ contains $\Qcyc$,
and we define

$$G=G(\Fin/\Q), \quad
H=G(\Fin/\Qcyc), \quad
\Ga=G(\Qcyc/\Q).
\eqno{(96)}$$
Thus $G$ is a compact $p$-adic Lie group
to which the abstract theory  developed in the
first four sections of this paper applies.
In fact, the structure of $G$ is well known.
Indeed, $G$ can be identified with
a closed subgroup of $GL_2(\Zp)=\Aut(\Epinf)$.
When $E$ admits complex multiplication,
$G$ has dimension 2, and when $E$ does not admit
complex multiplication, $G$ has
dimension 4 \cite{Se3}.\par

If $L$ is any algebraic extension of $\Q$,
we recall that the classical Selmer group
$\cS(E/L)$ is defined by

$$\cS(E/L)=\Ker(H^1(L, \Epinf) \to
\prod_w(H^1(L_w,E(\bar{L}_w))),$$
where $w$ runs over all non-archimedean places of $L$,
and $L_w$ denotes the union of the completions
at $w$ of all finite extensions of $\Q$ contained
in $L$. We write

$$X(E/L)=\Hom(\cS(E/L), \Qp/\Zp)$$
for the compact Pontrjagin dual of the
discrete abelian group $\cS(E/L)$.
We shall mainly be interested in the case in which
$L$ is Galois over $F$, in which case both
$\cS(E/L)$ and $X(E/L)$ have a natural left action
of $G(L/F)$, which extends to a left action of the
whole Iwasawa algebra $\Lambda(G(L/F))$.
It is easy to see that $X(E/L)$ is always a finitely
generated $\Lambda(G(L/F))$-module.\par

We impose for the rest of the paper the following
\bigskip
\newline
{\bf {Hypotheses on $p$.}}
\qquad $p \geq 5$ and $E$ has good ordinary
reduction at $p$.\par

\bigskip

The condition $p \geq 5$ guarantees that
$G$ has no element of order $p$,
since $G$ is a closed subgroup of $GL_2(\Zp)$.
It is a basic question as to when, assuming
our hypotheses on $p$,
the dual Selmer group $X(E/\Fin)$ belongs to the
category $\MSaG$, or equivalently is
$S^*$-torsion, where $S^*$ is the
Ore set defined in \S 3.\par

\begin{conj}
Assuming our hypotheses on $p$,
$X(E/\Fin)$ belongs to the category $\MSaG$.
\end{conj}

We now briefly discuss what is known at present about Conjecture 5.1.
Indeed, all we know is related to the following
much older conjecture of Mazur \cite{Ma}.

\begin{conj}
Assume $E$ has good ordinary reduction at $p$.
For each finite extension $F$ of $\Q$,
$X(E/\Fcyc)$ is $\La(\GaF)$-torsion.
\end{conj}

One might hope that Conjecture 5.1 is equivalent to
knowing Conjecture 5.2 for all
finite extensions $F$ of $\Q$
contained in $\Fin$.
At present,  we can only prove some
partial results in this direction,
which we now describe. If
$F$ is a finite extension of $\Q$
contained in $\Fin$,
we put

$$G_F=G(\Fin/F),
\quad
H_F=G(\Fin/\Fcyc).
\eqno{(98)}$$
Let $X$ be a finitely generated $\LaG$-module.
If $L$ is a finite extension of $\Q$ in
$\Fin$ such that $G_L$ is pro-$p$,
the $\mu$-invariant of $X$ viewed as
a $\La(G_L)$-module, is defined
in \cite{Ho} and \cite{Ve2}.
Similarly, if $Z$ is any finitely generated
$\La(\GaF)$-module, we write
$\mu_{\GaF}(Z)$ for its classical $\mu$-invariant.
Also, for any field $K$ in $\Fin$, we define
$$Y(E/K)=X(E/K)/X(E/K)(p).
\eqno{(99)}$$

\begin{lem}
Assume that $X(E/\Fin)$ belongs
to $\MSaG$.
Then, for each finite extension $F$ of
\:$\Q$
contained in $\Fin$,
we have
$\rm{(i)}$ $X(E/\Fcyc)$ is
$\La(\GaF)$-torsion$;$
$\rm{(ii)}$ $H_i(H_F, X(E/\Fin))=0$ for all
$i \geq 1;$
$\rm{(iii)}$ $H_i(H_F,Y(E/\Fin))$ is finite
for $i=1,2$, and
zero for $i > 2$.
Moreover, we have
$\mu_{G_L}(X(E/\Fin))
=\mu_{\Ga}(X(E/\Lcyc))$ for each
finite extension $L$ of $\Q$ in $\Fin$
such that $G_L$ is pro-$p$.
\end{lem}

\begin{pf}
We assume that $X(E/\Fin)$
belongs to $\MSaG$,
so that $Y(E/\Fin)$ is a finitely generated
$\LaH$-module.
Take $F$ to be any finite extension of
$\Q$ contained in $\Fin$.
Then $Y(E/\Fin)_{H_F}$ is clearly a
finitely generated $\Zp$-module,
because $Y(E/\Fin)$ is then a finitely
generated $\La(H_F)$-module.
Now we have the commutative diagram
with exact rows

$$
\begin{CD}
X(E/\Fin)_{H_F} @>>>Y(E/\Fin)_{H_F}@>>> 0\\
@VV{r_F}V @VV{s_F}V\\
X(E/\Fcyc)@>>>Y(E/\Fcyc)
@>>> 0,
\end{CD}
$$
where $r_F$ is the dual of the restriction map
$t_F$ from $\cS(E/\Fcyc)$ to $\cS(E/\Fin)^{H_F}$,
and $s_F$ is the corresponding induced map.
We claim that the cokernel of $r_F$ is finite.
Indeed, the cokernel of $r_F$ is dual to the kernel of
the restriction map $t_F$,
and $\Ker(t_F)$ is contained in $H^1(H_F, \Epinf)$.
But this latter group is proven to be finite in
\cite{CS1}. Since $\Coker(r_F)$ is finite,
it follows that $\Coker(s_F)$ is also finite.
As was remarked above, $Y(E/\Fin)_{H_F}$ is
a finitely generated $\Zp$-module,
and so we see that
$Y(E/\Fcyc)$
is a finitely generated $\Zp$-module.
Thus $\XEFcyc$ is $\La(\GaF)$-torsion,
proving (i).
To prove the remainder of Lemma 5.3,
we invoke the arguments of
\cite{CSS}, and do not repeat
the detailed proofs of these here.
Since $\XEFcyc$ is $\La(\GaF)$-torsion,
Lemma 2.5 of \cite{CSS} shows that
$H_1(H_F,\XEFin)=0$. As
we have shown $\XEFcyc$ is $\La(\Ga_F)$-torsion
for all $F$, Lemma 2.1 and Remark 2.6 of \cite{CSS}
prove that $H_i(H_F,X(E/\Fin))=0$
for all $i \geq 2$,
establishing (ii).
Turning to (iii),
an entirely similar argument to that given
in the proof of Proposition 2.13
of \cite{CSS} gives that $H_3(H_F,Y(E/\Fin))=0$,
and that $H_i(H_F,Y(E/\Fin))$ $(i=1,2)$ are
killed by a power of $p$.
Hence, as the $H_i(H_F,Y(E/\Fin))$ $(i \geq 0)$ are
finitely generated $\Zp$-modules,
we conclude that $H_i(H_F,Y(E/\Fin))$ must be
finite for $i=1,2$.
The final assertion of Lemma 5.3 now follows
from formula (25) of Proposition 2.13 of \cite{CSS}, on noting that
the $\mu$-invariant of the $\La(\Ga_L)$-module $H_0(H_L,Y(E/\Fin))$
is zero, because this module is finitely generated over $\Zp$.
This completes the proof of Lemma 5.3.
\end{pf}

\begin{lem}
\label{l5.4}
Assume that $\XEFcyc$ is $\La(\GaF)$-torsion
for all finite extensions $F$ of $\Q$ in $\Fin$.
Suppose, in addition, that there exists a finite
extension $L$ of $\Q$ contained in $\Fin$ satisfying$:-$
$\rm{(i)}$ $G_L$ is pro-$p$,
$\rm{(ii)}$ $\mu_{G_L}(\XEFin)=\mu_{\Ga_L}(\XELcyc)$,
and
$\rm{(iii)}$ $H_1(H_L,\YEFin)$ is finite.
Then $\XEFin$ belongs to $\MSaG$.
\end{lem}

\begin{pf}
We first observe that
$H_0(H_L,\YEFin)$ is a finitely generated torsion
$\La(\Ga_L)$-module, assuming the hypotheses of the
lemma. Indeed, $H_0(H_L, \YEFin)$ is a quotient of
$H_0(H_L, \XEFin)$, and we claim that this latter module is
a finitely generated torsion $\La(\Ga_L)$-module.
This is because we are assuming that the finitely $\La(\Ga_L)$-module
$\XELcyc$ is $\La(\Ga_L)$-torsion,
and it is shown in \cite{CH} that the kernel of the natural map

$$r_L:\XEFin_{H_L} \to \XELcyc,$$
which is the dual of the restriction map, is a finitely
generated $\Zp$-modules.
Hence $\Ker(r_L)$ is a finitely generated torsion
$\La(\Ga_L)$-module,
and our claim follows.
The delicate part of the proof is to now show that

$$\mu_{\Ga_L}(H_0(H_L,\YEFin))=0.
\eqno{(100)}
$$
Indeed, (100) implies immediately that
$\YEFin$ is a finitely generated
$\LaH$-module, completing
the proof of the lemma. This is because
(100) shows that $H_0(H_L, \YEFin)$ is
a finitely generated $\Zp$-module, and hence,
as $H_L$ is pro-$p$, Nakayama's lemma
gives that $\YEFin$ is finitely generated over
$\La(H_L)$. To prove (100), we invoke
the full force of Proposition 2.13 of
\cite{CSS}.
As remarked earlier, our assumption that
$\XEFcyc$ is $\La(\GaF)$-torsion for
all finite extensions $F$ of $\Q$ in $\Fin$ implies that all the
hypotheses of Proposition 2.13 of \cite{CSS} are valid.
Hence formula (25) of Proposition 2.13 of \cite{CSS} holds.
Inserting conditions (ii) and (iii) in formula (25) of \cite{CSS}, 
we conclude
that the $\mu_{\Ga_L}$-invariants of 
$H_0(H_L,\YEFin)$ and
$H_2(H_L,\YEFin)$
must add up to zero. Thus, as both are non-negative
integers, both of these $\mu$-invariants must be zero, proving (100) in
particular.
\end{pf}

\begin{cor}
Assume that $p \geq 5$, $E$ has good ordinary reduction at $p$,
and $E$ admits complex multiplication. Then $\XEFin$ belongs
to $\MSaG$ if and only if there exists a finite extension $L$
of $\Q$ in $\Fin$ such that $G_L$ is pro-$p$ and
$\mu_{\Ga_L}(\XELcyc)=0$.
\end{cor}

\begin{pf}
It is a well known consequence of the
proof of the two variable main conjecture for $E$ over
$\Fin$ (see \cite{Ru}) that $\XEFcyc$ is
$\La(\Ga_F)$-torsion for every finite extension $F$ of
$\Q$ in $\Fin$. Take $L$
to be any finite extension of $\Q$ contained in $\Fin$
such that $G_L$ is pro-$p$. Then it is well known that
the thesis of Schneps \cite{Sc} implies that $\mu_{G_L}(\XEFin)=0$.
This in turn implies that $\XEFin(p)=0$ because
$\XEFin$ has no non-zero $\La(G_L)$-pseudo-null
submodule by \cite{Pe}.
We then have $H_i(H_L,\YEFin)=0$ for all $i \geq 1$
because $\YEFin=\XEFin$. The assertion of the corollary is now
clear from Lemmas 5.3 and 5.4.
\end{pf}

The following is the one practical
criterion we know at present for proving in some
concrete examples that $\XEFin$ belongs to $\MSaG$.

\begin{prop}
Assume that there exists a finite
extension $L$ of $\Q$ contained in
$\Fin$, and an elliptic curve $E'$ defined over $L$,
as follows$:-$
$\rm{(i)}$ $G_L$ is pro-$p$,
$\rm{(ii)}$ $E'$ is isogenous to $E$ over $L$,
and
$\rm{(iii)}$ $\XEpLcyc$ is finitely generated
$\Zp$-module. Then $\XEFin$ belongs to $\MSaG$.
\end{prop}

It is very probable that the hypotheses of Proposition 5.6
are valid for most elliptic curves $E$ over $\Q$ and
most primes $p$ of good ordinary reduction.
However, in our present state of knowledge, condition
(iii) is difficult to verify in numerical examples.
Here is an example of one isogeny class of curves
to which it can be applied.\par

\bigskip

{\bf {Example.}}
Consider the three elliptic curves over
$\Q$ of conductor 11, which we denote
by $E_i$ $(i=0,1,2)$. Here $E_1$
is given by (70), and the other two are given by
the equations

$$E_0:y^2+y=x^3-x^2-10x-20$$

$$E_2:y^2+y=x^3-x^2-7820x-263580$$
All three curves are isogenous over
$\Q$, and so the field $\Fin$ is the same for the three
curves and all primes $p$.
Taking $p=5$, and $L=\Q(\mu_5)$, it is well known
\cite{Ru} that $G_L$ is pro-5 and $X(E_1/\Lcyc)=0$.
Hence Proposition 5.6 shows that
$X(E_i/\Fin)$ belongs to $\MSaG$ for $p=5$ and $i=0,1,2$.
At present, we do not know how to prove that $X(E_i/\Fin)$
belongs to $\MSaG$ for any good ordinary prime $p>5$.\par

\bigskip

We now prove Proposition 5.6. Pick an isogeny
$f:E' \to E$, which is defined over $L$.
Note that $\Fin=L(E'_{p^{\infty}})$.
The natural homomorphism

$$r_L:\XEpFin_{H_L} \to \XEpLcyc,$$
which is the dual of the restriction map,
has a kernel which is a finitely generated
$\Zp$-module (see \cite{CH}). Applying Nakayama's lemma,
we deduce from (i) and (iii) that
$\XEpFin$ is finitely generated over
$\La(H_L)$. On the other hand, it is easily seen that
our isogeny $f$ induces a $\La(G_L)$-homomorphism from
$\XEFin$ to $\XEpFin$, whose kernel is
killed by a power of $p$. Thus $\XEFin$ belongs
to $\MSaG$, and the proof of Proposition 5.6 is complete.
\qed\par

\bigskip

We end this paper by conjecturing the existence
of a $p$-adic $L$-function attached to $E$ over
$\Fin$, and formulating a main conjecture
which relates this $p$-adic $L$-function to
$\XEFin$. We refer
the reader to \cite{KF} for motivation
for the definition we have chosen for our $p$-adic $L$-function.
If $q$ is any prime number, we write
$\Frobq$ for the Frobenius automorphism of $q$ in
$\GalQq/I_q$, where, as usual, $I_q$ denotes the
inertia subgroup. To fix the definition of our local and
global $\eps$-factors, we follow the conventions and
normalizations of \cite{De1}. In particular, we choose
the sign of the local reciprocity map so that $q$ is
mapped to $\Frobq^{-1}$. We take the
Haar measure on the ad\`ele group of $\Q$ which is
the usual Haar measure on $\R$, and, for each
prime number $q$, gives $\Zq$ volume 1.
We also fix the unique complex character of the ad\`ele
group of $\Q$ whose infinite component is
$x \mapsto \exp(2\pi i x)$, and whose component at
a finite prime $q$ is $x \mapsto  \exp(-2 \pi i x)$. Now let
$\rho$ denote any Artin representation of $G$.
In our earlier $p$-adic theory, we viewed $\rho$ as being
realized over the ring of integers of some finite extension
of $\Qp$. However, since $\rho$ is
an Artin representation of $G$,
finite group theory tells us that there exists a finite
extension of $\Q$, which we denote
by $\Krho$, such that $\rho$ can be realized
in a finite dimensional
vector space $\Vrho$ over $\Krho$.
We first recall the definition of the complex Artin
$L$-function $L(\rho,s)$.
It is given by the Euler product

$$L(\rho,s)=\prod_q P_q(\rho,q^{-s})^{-1},$$
where $P_q(\rho,T)$ is the polynomial

$$P_q(\rho,T)=\det(1-\Frobq^{-1}.T | \Vrho^{I_q}).$$
For each prime number $q$, we write $\eqrho$ for the local
$\eps$-factor of $\rho$ at $q$,
normalized as in \cite{De1}.\par

The complex $L$-function which is obtained by
twisting $E$ by the Artin representation $\rho$
is defined as follows. For each prime number $l$,
let

$$T_l(E)=\lim_{\longleftarrow}E_{l^n},\quad
V_l(E)=T_l(E)\otimes_{\Zl}\Ql, \quad
H^1_l(E)=\Hom(\Vl(E), \Ql).
\eqno{(101)}$$
Moreover, we fix some prime $\la$ of $\Krho$ above $l$,
and put $\Vrhola=\Vrho \otimes_K K_{\rho,\la}$.
Then

$$L(E,\rho,s)=\prod_qP_q(E,\rho,q^{-s})^{-1},
\eqno{(102)}$$
where $P_q(E,\rho,T)$ is the polynomial

$$P_q(E,\rho,T)=
\det(1-\Frobq^{-1}.T|(
H^1_l(E)\otimes_{\Ql}\Vrhola)^{I_q});$$
here $l$ is any prime number different from $q$.
The Euler product $L(E,\rho,s)$ converges only for
$\Re(s)>3/2$, and the only thing known about its
analytic continuation at present is that it has a meromorphic
continuation when $\rho$ factors through a soluble extension
of $\Q$. We will assume the analytic continuation
of $L(E, \rho,s)$ to $s=1$ for all Artin characters
$\rho$ of $G$ in what follows. The point $s=1$ is critical
for $L(E,\rho,s)$ for all Artin $\rho$,
in the sense of \cite{De2} and the period
conjecture of \cite{De2} asserts the following in
this case. Fix a global minimal Weierstrass equation
for $E$ over $\Z$, and let $\om$ denote the N\'eron differential
of this equation. Let $\ga^+$ (resp. $\ga^{-}$) denote a
generator for the subspace of $H_1(E(\C), \Z)$ fixed by
complex conjugation (resp. the subspace on which complex
conjugation acts by -1). We define

$$\Om_+(E)=\int_{\ga^+}\om, \quad
\Om_-(E)=\int_{\ga^{-}}\om.$$
Moreover, let $d^+(\rho)$ denote the dimension of the
subspace of $\Vrho$ on which complex conjugation acts
by $+1$, and $d^-(\rho)$ the dimension of the
subspace on which complex conjugation acts by -1.
According to the period conjecture of \cite{De2},
we have

$$\frac{L(E,\rho,1)}{\Om_+(E)^{d^+(\rho)}\Om_-(E)^{d^-(\rho)}}
\in \barQ
\eqno{(103)}$$
for all Artin representation $\rho$ of $G$. We shall
assume (103) to define our conjectural $p$-adic $L$-function.
We also tacitly suppose we have fixed embeddings of $\barQ$
into both $\C$ and $\barQp$.\par

Let $j_E$ denote the $j$-invariant of $E$, and let
$R$ denote the set consisting of the prime $p$
and all prime numbers $q$ with
$\ord_q(j_E)<0$. We define

$$L_R(E,\rho,s)=\prod_{q \notin R}P_q(E, \rho, q^{-s})^{-1}.
\eqno{(104)}$$
Finally, we put

$$p^{f_{\rho}}={\text {$p$-part of the conductor of $\rho$.}}
\eqno{(105)}$$
Also, since $E$ is ordinary at $p$, we have

$$1-a_pX+pX^2=(1-uX)(1-wX), \quad
u \in \Zp^{\times};
\eqno{(106)}$$
here, as usual,
$p+1-a_p=\sharp(\tilde{E}_p(\F_p))$, where
$\tilde{E}_p$ denotes the reduction of $E$ modulo $p$.
We recall that $\hat{\rho}$ denotes
the contragredient representation of $\rho$.\par

\begin{conj}
Assume that $p \geq 5$ and that $E$ has good
ordinary reduction at $p$. Then there
exists $\cLE$ in
$K_1(\LaGSa)$ such that,
for all Artin representations $\rho$ of $G$,
we have $\cLE(\rho) \ne \infty$, and

$$\cLE(\rho)=\frac{L_R(E,\rho,1)}
{\Om_+(E)^{d^+(\rho)}\Om_-(E)^{d^-(\rho)}}\cdot
\eprho \cdot
\frac{P_p(\hat{\rho},u^{-1})}
{P_p(\rho,w^{-1})} \cdot
u^{-f_{\rho}};
\eqno{(107)}$$
here $\eprho$ denotes the local $\eps$-factors at $p$ attached
to $\rho$.
\end{conj}

We remark that, when $E$ admits complex multiplication,
Conjecture 5.7 is true, and
can be deduced from the existence of the two variable
$p$-adic $L$-function of
Manin-Vishik, Katz, and Yager
attached to $E$ at the ordinary prime $p$.
When $E$ does not admit complex multiplication,
the only evidence we have at present to support
Conjecture 5.7 is some interesting but still
fragmentary numerical data due to Balister \cite{Bal}
(unpublished) and T.Dokchitser-V.Dokchitser \cite{27}.

\begin{conj}[The main conjecture]
Assume that $p \geq 5$, $E$ has good ordinary reduction
at $p$, and $\XEFin$ belongs to
$\MSaG$. Granted Conjecture $5.7$, the $p$-adic
$L$-function $\cLE$ in $K_1(\LaGSa)$ is a characteristic
element of $\XEFin$.
\end{conj}

It is plain that such a main conjecture will have many
deep consequences for the arithmetic of $E$ over
$\Fin$ and its subfields.
We do not enter into a discussion of these here, beyond
mentioning the following almost immediate corollaries.\par

\begin{cor}
Assume Conjecture $5.8$. Then, for each
Artin representation $\rho$ of $G$,
$\chi(G, \twrho(\XEFin))$ is finite if and only if
$L(E, \hat{\rho},1) \ne 0$.
\end{cor}

\begin{cor}
Assume Conjecture $5.8$. Let $\rho$ be an Artin
representation of $G$ such that $L(E, \hat{\rho},1) \ne 0$.
As in Theorem $3.6$, let $\mrho$ be the degree of the
quotient field of $O$ over $\Qp$, where
$\rho$ is given by $(12)$.
Then $\chi(G, \twhatrho(M))$ is equal to the
$\mrho$-th power of the inverse of the $p$-adic valuation of
the right hand side of $(107)$.
\end{cor}

We make the following remark about the integrality
properties of $\cLE$ and any characteristic element
$\xi_E$ of
$\XEFin$. On the one hand,
the conjectured holomorphy of $L(E,\rho,s)$ at
$s=1$ gives, via Conjecture 5.7, the assertion that
$\cLE(\rho)\ne \infty$ for every Artin representation
$\rho$ of $G$. On the other hand,
thanks to Lemma 5.3 and Theorem 3.8, the conjecture
that $\XEFin$ belongs to $\MSaG$ implies that
$\xi_E(\rho)\ne \infty$ for every Artin representation
$\rho$ of $G$. Note that Case 1 of Conjecture 4.8 would
then imply that both $\cLE$ and $\xi_E$ belong to
the image of the natural map

$$\LaG[\frac{1}{p}]^{\times} \cap \LaGSa
\hookrightarrow \LaGSa \to
K_1(\LaGSa).$$
To prove Corollary 5.9, we observe that when
$\XEFin$ belongs to $\MSaG$, Lemma 5.3 and Theorem 3.8
show that $\xi_E(\rho)\ne 0$ if and only if
$\chi(G, \twhatrho(\XEFin))$ is finite.
If we assume Conjecture 5.8, we can take
$\xi_E=\cLE$, and then it is clear from
(107) that $\cLE(\rho) \ne 0$ if and only if
$L(E,\rho,1)\ne 0$. Corollary 5.10 is then an
immediate consequence of (107) and formula (36) of
Theorem 3.6.
\qed\par

\bigskip

Needless to say, the evidence in favour
of Conjecture 5.8 is still very limited. When $E$ admits complex
multiplication, Conjecture 5.8 is true
provided we assume that $\XEFin$ belongs to $\MSaG$
(see Corollary 5.5). Indeed
it can be deduced from the proof by
Yager \cite{28} and Rubin \cite{Ru}
of what is usually called the two variable main conjecture
(recall that in this case $G$ is a $p$-adic Lie group of
dimension 2).
When $E$ does not admit complex multiplication, the remarkable
calculations of $L(E,\rho, 1)$ made in \cite{27}
for certain Artin representations of small degree and
small primes $p$ provide some numerical evidence in
favour of Conjecture 5.8.
We end by giving their results for $E=X_1(11)$, $p=5$,
and the two Artin representations $\rho_1$ and $\rho_2$
of degree 4 which are discussed at the end of \S 3.\par

\bigskip

{\bf {Example.}} We use the same notation as at the end of
\S 3. The following data is calculated in \cite{27}.
We have

$$d^+(\rho_i)=d^-(\rho_i)=2
\quad(i=1,2)
\eqno{(108)}$$

$$f_{\rho_1}=3, \quad
f_{\rho_2}=5
\eqno{(109)}$$

$$P_5(\rho_1,X)=1-X, \quad
P_5(\rho_2,X)=1
\eqno{(110)}$$

$$P_5(E,\rho_1,X)=1-X+5X^2, \quad
P_5(E,\rho_2,X)=1
\eqno{(111)}$$

$$P_{11}(E,\rho_1,X)=1, \quad
P_{11}(E,\rho_2,X)=1
\eqno{(112)}$$

$$|e_5(\rho_1)|_5=5^{-3/2}, \quad
|e_5(\rho_2)|_5=5^{-5/2}
\eqno{(113)}$$
Moreover, $\sharp(\tilde{E}(\F_5))=5$, and so

$$1-a_5X+5X^2=1-X+5X^2
\eqno{(114)}$$
Finally, even though
$P_q(E,\rho_i,X)$
$(i=1,2)$ has degree 8 for all $q\ne 5,11$, the remarkable
calculations of \cite{27} show that

$$\frac{L(E,\rho_1,1)}
{(\Om_+(E)\Om_-(E))^2}=
\frac{-2^2}{{11}^2\sqrt{5}},
\quad
\frac{L(E,\rho_2,1)}
{(\Om_+(E)\Om_-(E))^2}=
\frac{-2^6}{{11}^25\sqrt{5}}.
\eqno{(115)}$$
Thus, since both of these $L$-values are non-zero,
Proposition 3.11 shows that Corollary 5.9
holds for $\rho_1$ and $\rho_2$. Finally,
one can use
the above data to calculate the right hand side of
(107), up to a 5-adic unit, for
$\rho_1$ and $\rho_2$. Using (74) of Proposition 3.11,
and recalling that $\rho_1$ and $\rho_2$ can
both be realized over $\Z_5$, we deduce that
Corollary 5.10 is also valid for $\rho_1$
and $\rho_2$.\par

\bigskip

\begin{tabbing}
xxxx\=
xxxxxxxxxxxxxxxxxxxxxxxxxxxxxxx\= xxxxxxxxxxxxxxxxxxxxxxxxxxxxxxxxxxxxxx\=\kill
\> John Coates \>                       Takako Fukaya\\
\> DPMMS,\>                             Faculty of Business and Commerce,\\
\> University of Cambridge,\>           Keio University,\\
\> Centre for Mathematical Sciences,\>  Hiyoshi, Kohoku-ku,\\
\> Wilberforce Road, \>                 Yokohama, 223-8521, Japan\\
\> Cambridge CB3 0WB, England \>        takakof@hc.cc.keio.ac.jp\\
\> J.H.Coates@dpmms.cam.ac.uk\>         \\               

\end{tabbing}

\begin{tabbing}

xxxx\=
xxxxxxxxxxxxxxxxxxxxxxxxxxxxxxx\= xxxxxxxxxxxxxxxxxxxxxxxxxxxxxxxxxxxxxx\kill
\> Kazuya Kato\>                     Ramdorai Sujatha\\
\> Department of Mathematics,\>      School of Mathematics,\\
\> Kyoto University,\>               Tata Institute of Fundamental Research,\\
\> Kitashirakawa,\>                  Homi Bhabha Road,\\
\> Kyoto 606-8502, Japan\>           Mumbai 400 005, India\\
\> kazuya@math.kyoto-u.ac.jp\>       sujatha@math.tifr.res.in\\
\end{tabbing}

\begin{tabbing}
xxxx\=
xxxxxxxxxxxxxxxxxxxxxxxxxxxxxxx\= xxxxxxxxxxxxxxxxxxxxxxxxxxxxxxxxxxxxxx\kill
\> Otmar Venjakob\\
\> Universit\"at Heidelberg,\\
\> Mathematisches Institut,\\
\> Im Neuenheimer Feld 288,\\
\> D-69120 Heidelberg, Germany\\
\>otmar@mathi.uni-heidelberg.de
\end{tabbing}


\medskip





\begin{thebibliography}{99}

\bibitem{Ba}{\sc Bass, H.},
\newblock
{\em Algebraic $K$-theory},
Benjamin, New York (1968).

\bibitem{Bal}{\sc Balister, P.},
\newblock
{\em Congruences between special values of 
$L$-functions} (unpublished).



\bibitem{BCA}{\sc Bourbaki, N.},
\newblock
{\em Commutative Algebra},
Hermann, Paris (1965).



\bibitem{Bru}{\sc Brumer, A.},
\newblock
{\em Pseudocompact algebras, profinite groups and class
formations},
J. of Algebra {\bf 4} (1966) 442-470.

\bibitem{CH}
{\sc Coates, J.,
{\rm and} Howson, S.},
\newblock
{\em Euler characteristics and elliptic curves
{\rm II}},
J. Math. Soc. Japan Proc. {\bf 53} (2001) 175--235.

\bibitem{CSS}{\sc Coates, J., Schneider, P.,
{\rm and} Sujatha, R.},
\newblock
{\em Links between cyclotomic and $GL_2$ Iwasawa theory},
Doc. Math. Extra Volume : Kazuya Kato's 50th birthday (2003) 187--215.

\bibitem{31} 
{\sc Coates, J., Schneider, P.,
{\rm and} Sujatha, R.},
\newblock
{\em Modules over Iwasawa algebras},
J. Inst. Math. Jussieu {\bf 2} (2003)  no. 1, 73--108.


\bibitem{CS1}
{\sc Coates, J.,
{\rm and} Sujatha, R.},
\newblock
{\em Euler-Poincar\'e characteristics of abelian
varieties},
CRAS {\bf 329} S\'erie I (1999) 309--313.

\bibitem{CS}
{\sc Coates, J.,
{\rm and} Sujatha, R.},
\newblock
{\em Galois cohomology of elliptic curves},
TIFR Lecture notes series,
Narosa Publishing House (2000).


\bibitem{De1}
{\sc Deligne, P.},
\newblock
{\em Les constantes des \'equations fonctionnelles des fonctions
$L$},
Modular functions of one variable II LNM
{\bf 349} Springer (1973) 501--597.

\bibitem{De2}
{\sc Deligne, P.},
\newblock
{\em Valeurs de fonctions $L$ et p\'eriodes d'int\'egrales},
Proc. Sympos. Pure Math., XXXIII
Automorphic forms, representations and $L$-functions,
Part 2 Amer. Math. Soc. (1979) 313--346.


\bibitem{27}
{\sc Dokchitser, T., {\rm and} Dokchitser V.},
\newblock
Paper in preparation.


\bibitem{Fi}
{\sc Fisher, T.},
\newblock
{\em Descent calculations for the elliptic curves
of conductor $11$},
Proc. London Math. Soc. {\bf 86} (2003) 583--606.


\bibitem{KF}
{\sc Fukaya, T., {\rm and}
Kato, K.},
\newblock
{\em A formulation of conjectures on
$p$-adic zeta functions in non-commutative
Iwasawa theory}, preprint (2003).


\bibitem{Gr}
{\sc Greenberg, R.},
\newblock
{\em On the structure of certain Galois groups},
Invent. Math. {\bf 47} (1978) 85--99.


\bibitem{Ho}
{\sc Howson, S.},
\newblock
{\em Euler characteristics as invariants
of Iwasawa modules},
Proc. London Math. Soc. {\bf 85} (2002) 634--658.




\bibitem{La}
{\sc Lazard, M.},
\newblock
{\em Groupes analytiques $p$-adiques},
Publ. Math. IHES {\bf 26} (1965) 389--603.

\bibitem{Ma}
{\sc Mazur, B.},
\newblock
{\em Rational points of abelian varieties
in towers of number fields},
Invent. Math. {\bf 18} (1972) 183--266.


\bibitem{MR}
{\sc McConnell, ~J. ~C., {\rm {and}} Robson, ~J.~C.},
\newblock
{\em Noncommutative Noetherian Rings},
Graduate Studies in Math. {\bf 30} AMS (1987).

\bibitem{NSW} 
{\sc Neukirch, J., Schmidt, A., {\rm and} 
Wingberg, K.},
\newblock 
{\em Cohomology of number fields},  
Grundlehren der Mathematischen Wissenschaften 
{\bf 323}
Springer (2000).


\bibitem{Pe}
{\sc Perrin-Riou, ~B.},
\newblock
{\em Groupes de Selmer d'une courbe elliptique \`a multiplication
complexe},
Comp. Math. {\bf 43} (1981) 387--417.


\bibitem{Ru}
{\sc Rubin, ~K.},
\newblock
{\em The ``main conjectures" of Iwasawa theory for imaginary
quadratic fields},
Invent. Math. {\bf 103} (1991) 25--68.


\bibitem{Sc}
{\sc Schneps, ~L.},
\newblock
{\em On the $\mu$-invariant of $p$-adic $L$-functions},
J. Number Theory {\bf 25} (1987) 20-33.




\bibitem{Se2}
{\sc Serre, J.-P.},
\newblock
{\em Sur la dimension cohomologique de groupes profinis},
Topology {\bf 3} (1965) 413--420, Oeuvres II 264-271.



\bibitem{Se3}
{\sc Serre, J.-P.},
\newblock
{\em Properi\'et\`es Galoisiennes
des points d'ordre fini des courbes elliptiques},
Invent. Math. {\bf 15} (1972) 259--331.

\bibitem{30}
{\sc Serre, J.-P.}, 
\newblock
{\em Alg\'ebre Locale, Multipicit\'es}, 
3rd ed. LNM {\bf 11}
Springer (1975).


\bibitem{Se}
{\sc Serre, J.-P.},
\newblock
{\em Linear representations of finite groups},
Graduate Texts in Mathematics {\bf 42}
Springer (1977).


\bibitem{Sw}
{\sc Swan, R.},
\newblock
{\em Algebraic $K$-theory}, LNM {\bf 76} Springer (1968).

\bibitem{Va}
{\sc Vaserstein, L. N.},
\newblock
{\em On stabilization for general linear groups over a
ring}, 
Math. USSR Sbornik {\bf 8} (1969) 383-400.

\bibitem{Va2}
{\sc Vaserstein, L. N.},
\newblock
{\em On the Whitehead Determinant for Semi-local Rings},
Preprint (2004).




\bibitem{Ve2}
{\sc Venjakob, O.},
\newblock
{\em On the structure theory of the Iwasawa
algebra of a $p$-adic Lie group},
J. European Math. Soc. {\bf 4} (2002) 271--311.


\bibitem{32}
{\sc Venjakob, O.} (with an appendix by D. Vogel), 
\newblock
{\em A non-commutative Weierstrass preparation theorem and applications to
Iwasawa theory}, 
J. Reine Angew. Math. {\bf 559} (2003) 153--191.


\bibitem{Ve1}
{\sc Venjakob, O.},
\newblock
{\em Characteristic elements in non-commutative Iwasawa theory},
Habilitationschschrift, Heidelberg University (2003).


\bibitem{28}
{\sc Yager, R.I.},
\newblock
{\em On two variable $p$-adic $L$-functions},
Ann. of Math. {\bf 115} (1982) 411--449.













\end{thebibliography}
\end{document}